\documentclass{amsart}
\usepackage{amssymb}
\newtheorem{thm}{Theorem}[section]
\newtheorem{proposition}[thm]{Proposition}

\newtheorem{lemma}[thm]{Lemma}

\title[Grothendieck  rings of basic classical  Lie superalgebras]
{Grothendieck  rings of basic classical  Lie superalgebras}
\author{A.N. Sergeev}
\address{Department of Mathematical Sciences,
Loughborough University, Loughborough LE11 3TU, UK and Department of Mathematics, Saratov State University, Astrakhanskaya 83, Saratov, 410012, Russia}
\email{A.N.Sergeev@lboro.ac.uk}
\author{A.P. Veselov}
\address{Department of Mathematical Sciences,
Loughborough University, Loughborough LE11 3TU, UK  and Moscow State University, Moscow 119899, Russia}
\email{A.P.Veselov@lboro.ac.uk}
\begin{document}
\maketitle

\begin{abstract}
The Grothendieck rings of finite dimensional representations of the basic classical Lie superalgebras are explicitly described in terms of the corresponding generalized root systems. 
We show that they can be interpreted as the subrings in the weight group rings invariant under the action of certain groupoids called super Weyl groupoids.
\end{abstract}

\tableofcontents

\section{Introduction}
The classification of finite-dimensional representations of the semisimple complex Lie algebras and related Lie groups is one of the most beautiful pieces of mathematics. In his essay \cite{Atiyah} Michael Atiyah mentioned representation theory of Lie groups as an example of a theory, which "can be admired because of its importance, the breadth of its applications, as well as its rational coherence and naturality." This classical theory goes back to fundamental work by \`Elie Cartan and Hermann Weyl and is very well presented in various books, of which we would like to mention the famous Serre's lectures \cite{Serre} and a nicely written Fulton-Harris course \cite{FH}. One of its main results can be formulated as follows (see e.g. \cite{FH}, Theorem 23.24):

%\medskip

 {\it The representation ring $R(\mathfrak{g})$ of a complex semisimple Lie algebra $\mathfrak{g}$ is isomorphic to the ring ${\mathbb{Z}}[P]^{W}$ of $W$-invariants in the integral group ring ${\mathbb{Z}}[P],$
 where $P$ is the corresponding weight lattice and $W$ is the Weyl group. The isomorphism is given by the character map
$Ch: R(\mathfrak{g}) \rightarrow {\mathbb{Z}}[P]^{W}.$}

The main purpose of the present work is to generalize this result to the case of basic classical complex Lie superalgebras. The class of basic classical Lie superalgebras was introduced by Victor Kac in his fundamental work \cite{Kac, Kac2}, where the basics of the representation theory of these Lie superalgebras had been also developed. The problem of finding the characters of the finite-dimensional irreducible representations turned out to be very difficult and still not completely resolved (see the important papers by Serganova \cite{Serga1, Serga2} and Brundan \cite{Brun} and references therein).
Our results may shed some light on these issues.

Recall that a complex simple Lie superalgebra $\mathfrak {g} = \mathfrak{g}_{0} \oplus \mathfrak{g}_{1}$ is called the {\it basic classical} if it admits a non-degenerate invariant (even) bilinear form and the representation of the Lie algebra $\mathfrak{g}_{0}$  on the odd part $ \mathfrak{g}_{1}$ is completely reducible. The class of these Lie superalgebras can be considered as a natural analogue of the ordinary simple Lie algebras. In particular, they can be described (with the exception of $A(1,1) = \mathfrak{psl}(2,2)$) in terms of Cartan matrix and generalized root systems (see \cite{Kac2, Serga}).

Let $\mathfrak {g}$ be such Lie superalgebra
different from $A(1,1)$ and $\mathfrak{h}$ be its Cartan subalgebra (which in this case is also Cartan subalgebra of the Lie algebra $\mathfrak{g}_{0}).$ Let $P_0 \subset \mathfrak{h}^*$ be the abelian group of weights of $\mathfrak{g}_{0}$, $W_0$ be the  Weyl group of $\mathfrak{g}_{0}$ and  ${\mathbb{Z}}[P_0]^{W_0}$ be the ring of $W_0$-invariants in the integral group ring ${\mathbb{Z}}[P_0].$
The decomposition of $\mathfrak {g}$ with respect to the adjoint action of $\mathfrak{h}$ gives the (generalized) root system $R$ of Lie superalgebra $\mathfrak {g}.$ By definition $\mathfrak {g}$ has a natural non-degenerate bilinear form on $\mathfrak{h}$ and hence on $\mathfrak{h}^*,$ which will be denoted as $( , ).$ In contrast to the theory of semisimple Lie algebras some of the roots $\alpha \in R$ are isotropic: $(\alpha, \alpha) = 0.$  For isotropic roots one can not define the usual reflection, which explains the difficulty with the notion of Weyl group for Lie superalgebras. A geometric description of the corresponding generalized root systems were found in this case by Serganova \cite{Serga}.

Define the following {\it ring of exponential super-invariants} $J(\mathfrak {g})$, replacing the algebra of Weyl group invariants in the classical case of Lie algebras:
\begin{equation}
\label{Jdefin}
J(\mathfrak{g})=\{ f\in\mathbb{Z}[ P_0]^{W_0}:  \, D_{\alpha}f\in (e^{\alpha}-1) \quad \text{for any isotropic root } \alpha \}
\end{equation}
where $(e^{\alpha}-1)$ denotes the principal ideal in $ \mathbb{Z}[ P_0]^{W_0}$ generated by $e^{\alpha}-1$ and the derivative $D_\alpha$ is defined by the property $D_\alpha (e^{\beta}) = (\alpha, \beta) e^{\beta}.$
This ring is a variation of the algebra of invariant polynomials investigated
for Lie superalgebras in \cite{Ber}, \cite{Kac3}, \cite{Serge0,Serge1}.
For the special case of the Lie superalgebra $A(1,1)$ one should slightly modify the definition because of the multiplicity 2 of the isotropic roots (see section 8 below).

Our main result is the following

%\medskip

{\bf Theorem.}  {\it The Grothendieck ring $K(\mathfrak{g})$ of finite dimensional representations of a basic classical Lie superalgebra $\mathfrak{g}$ is isomorphic to the ring  $J(\mathfrak {g}).$ The isomorphism is given by the supercharacter map}
$Sch: K(\mathfrak{g}) \rightarrow J(\mathfrak {g}).$

The fact that the supercharacters belong to the ring $J(\mathfrak {g})$ is relatively simple, but the proof of surjectivity of the supercharacter map is much more involved and based on classical Kac's results \cite{Kac, Kac2}. 

The elements of $J(\mathfrak {g})$ can be described as the invariants in the weight group rings under the action of the following groupoid $\mathfrak{W}$, which we call {\it super Weyl groupoid}. It is defined as a disjoint union $$\mathfrak{W}(R) = W_0 \coprod W_0 \ltimes \mathfrak T_{iso},$$ where $\mathfrak T_{iso}$ is the groupoid, whose base is the set $R_{iso}$ of all isotropic roots of $\mathfrak {g}$ and  the set of morphisms from $\alpha \rightarrow \beta$ with $\beta \neq \alpha$ is non-empty if and only if $\beta = - \alpha$ in which case it consists of just one element $\tau_{\alpha}$. This notion was motivated by our work on deformed Calogero-Moser systems \cite{SV}.

The group $W_0$ is acting on $\mathfrak T_{iso}$ in a natural way and thus defines a semi-direct product groupoid $W_0 \ltimes \mathfrak T_{iso}$ (see details in section 9). One can define a natural action of $\mathfrak{W}$ on $\mathfrak h$ with $\tau_{\alpha}$ acting as a shift by $\alpha$ in the hyperplane $(\alpha,x)=0.$ If we exclude the special case of $A(1,1)$ our Theorem can now be reformulated as the following version of the classical case:

 {\it The Grothendieck ring $K(\mathfrak{g})$ of finite dimensional representations of a basic classical Lie superalgebra $\mathfrak{g}$ is isomorphic to the ring  ${\mathbb{Z}}[P_0]^{\mathfrak{W}}$ of the invariants of the super Weyl groupoid $\mathfrak{W}.$}

An explicit description of the corresponding rings $J(\mathfrak {g})$ (and thus the Grothendieck rings) for each type of basic classical Lie superalgebra is given in sections 7 and 8. 
For classical series we describe also the subrings, which are the Grothendieck rings of the corresponding natural algebraic supergroups.

\section{Grothendieck rings of Lie superalgebras}

All the algebras and modules in this paper will be considered over the field of complex numbers $\mathbb C.$

Recall that {\it superalgebra} (or {\it  $\mathbb Z_{2}$-graded algebra})  is an associative algebra $A$ with a decomposition into direct sum $A=A_{0}\oplus A_{1}$, such that if $a\in A_{i}$ and $b\in A_{j}$ then $ab\in A_{i+j}$ for all $i,j\in \mathbb Z_{2}$.  We will write $p(a)=i\in \mathbb Z_{2}$ if $a\in A_{i}$.

A {\it module over superalgebra} $A$ is a vector space $V$ with a decomposition
 $V=V_{0}\oplus V_{1}$, such that if $a\in A_{i}$ and $v\in V_{j}$ then $av\in V_{i+j}$ for all $i,j\in \mathbb Z_{2}$. Morphism  of $A$-modules $f : V\rightarrow U$ is module homomorphism preserving their gradings:  $f(V_{i})\subset U_{i},\,i\in \mathbb Z_{2}$. 
 
 We have the parity change functor  $V\longrightarrow \Pi(V)$,
where $\Pi(V)_{0}=V_{1},\: \Pi(V)_{1}=V_{0}$, with the $A$ action $a*v=(-1)^{p( a)}av$.
If $A,B$ are superalgebras then $A\otimes B$ is a superalgebra with the multiplication
$$
(a_{1}\otimes b_{1})(a_{2}\otimes b_{2})=(-1)^{p( b_{1})p( a_{2})}a_{1}a_{2}\otimes b_{1} b_{2}
$$
The tensor product of $A$-module $V$ and  $B$-module $U$ is $A\otimes B$-module  $V\otimes U$   and
$$
(V\otimes U)_{0}=(V_{0}\otimes U_{0})\oplus(V_{1}\otimes U_{1}),\quad(V\otimes U)_{1}=(V_{1}\otimes U_{0})\oplus(V_{0}\otimes U_{1})
$$
with the action $a\otimes b(v\otimes u)=(-1)^{ b v}av\otimes bu$.

The  {\it Grothendieck group} $K(A)$ is defined (cf. Serre \cite{Serre2}) as the quotient of the free abelian group with generators given by  all isomorphism classes of finite dimensional $\mathbb Z_{2}$-graded $A$-modules by the subgroup generated by
$
[V_{1}]-[V]+[V_{2}]
$
for all exact sequences
$$
0\longrightarrow V_{1}\longrightarrow V \longrightarrow V_{2} \longrightarrow 0
$$
and by
$
[V]+[\Pi(V)]
$
for all $A$-modules $V.$

It is easy to see that the Grothendieck group $K(A)$ is a free $\mathbb Z$-module with the basis corresponding to the classes of the irreducible modules.

Let now $A = U(\mathfrak g)$ be the universal enveloping algebra of a Lie superalgebra $\mathfrak g$ (see e.g. \cite{Kac}) and $K(A)$ be the corresponding Grothendieck group. Consider the map
$$
\mathfrak g\longrightarrow U(\mathfrak g)\otimes U(\mathfrak g),\quad x\rightarrow x\otimes 1+1\otimes x
$$
One can check that this map is a homomorphism of Lie superalgebras, where on the right hand side we consider the standard Lie superalgebra structure defined for any associative algebra $A$ by the formula
$$[a,\,b] = ab - (-1)^{p(a)p(b)} ba.$$
Therefore one can define
for any two $\mathfrak g $-modules $V$ and  $U$ the $\mathfrak g $-module structure on $V\otimes U$. Using this we define the product  on $K(A)$ by the formula
$$[U] [V] = [U \otimes V].$$ Since all modules are finite dimensional this multiplication is well-defined on the Grothendieck group $K(A)$ and introduces the ring structure on it. The corresponding ring is called {\it Grothendieck ring of Lie superalgebra} $\mathfrak g$ and will be denoted $K(\mathfrak g).$

\section{Basic classical Lie superalgebras and generalized root systems}

Following Kac \cite{Kac, Kac2} we call Lie superalgebra $\mathfrak {g} = \mathfrak{g}_{0} \oplus \mathfrak{g}_{1}$ {\it basic classical} if

a) $\mathfrak {g}$ is simple,

b) Lie algebra $\mathfrak{g}_{0}$ is a reductive subalgebra of $\mathfrak {g},$

c) there exists a non-degenerate invariant even bilinear form on $\mathfrak {g}$.

 Kac proved that the complete list of basic classical Lie superalgebras, which are not Lie algebras, consists of Lie superalgebras of the type
$$A(m,n), \, B(m,n), \, C(n), \, D(m,n), \, F(4), \, G(3),\, D(2,1,\alpha).$$

In full analogy with the case of simple Lie algebras one can consider the decomposition of $\mathfrak {g}$ with respect to adjoint action of Cartan subalgebra $\mathfrak h:$
$$ \mathfrak {g} = \mathfrak {h} \oplus \left(\oplus \mathfrak {g}_{\alpha} \right),$$
where the sum is taken over the set $R$ of non-zero linear forms on $\mathfrak {h}$, which are called {\it roots} of $\mathfrak {g.}$ With the exception of the Lie superalgebra of type $A(1,1)$ the corresponding root subspaces $\mathfrak {g}_{\alpha}$ have dimension 1 (for $A(1,1)$ type the root subspaces corresponding to the isotropic roots have dimension 2).

It turned out that the corresponding root systems admit the following simple geometric description found by Serganova \cite{Serga}.

Let $V$ be a finite dimensional complex vector space with a
non-degenerate bilinear form $( , )$.

{\bf Definition} \cite{Serga}. The finite set $R\subset V\setminus\{0\}$ is
called a {\it generalized root system} if the following conditions
are fulfilled :

1) $R$ spans $V$ and $R=-R$ ;

2) if $\alpha,\beta\in R$ and  $(\alpha ,\alpha )\ne 0$ then
$\frac{2(\alpha ,\beta )}{(\alpha ,\alpha )}\in {\mathbb{Z}}$ and
$s_{\alpha}(\beta)=\beta -\frac{2(\alpha ,\beta )}{(\alpha ,\alpha
)}\alpha\in R$;

3) if $\alpha\in R$ and $(\alpha ,\alpha )=0$ then
there exists an invertible mapping $r_{\alpha} : R \rightarrow R$ such that
$r_{\alpha} (\beta) = \beta$ if $(\beta, \alpha) = 0$ and $r_{\alpha} (\beta) \in \{ \beta + \alpha, \, \beta - \alpha\}$ otherwise.

The roots $\alpha$ such that $(\alpha, \alpha) = 0$ are called
{\it isotropic.}  A generalized root system $R$ is called {\it reducible} if it can
be represented as a direct orthogonal sum of two non-empty generalized root
systems $R_{1}$ and
 $R_{2}$: $V=V_{1}\oplus V_{2}$, $R_{1}\subset V_{1}$, $R_{2}\subset V_{2}$,
$R=R_{1}\cup R_{2}$. Otherwise the system is called {\it irreducible.}

Any generalized root system has a partial symmetry described by
the finite group $W_0$ generated by reflections with respect
to the non-isotropic roots.

A remarkable fact proved by Serganova \cite{Serga} is
that classification list for the irreducible generalized root
systems with isotropic roots coincides with the root systems of the basic classical Lie superalgebras from the Kac list
(with the exception of $A(1,1))$ and $B(0,n)$). Note that the superalgebra $B(0,n)$ has no isotropic roots: its root system is the usual non-reduced system of $BC(n)$ type.

{\bf Remark.}  Serganova considered also a slightly wider notion \footnote {Johan van de Leur communicated to us that a similar notion was considered earlier by T. Springer, but his classification results were not complete \cite{Leur}.} of generalized root systems, when the property 3) is replaced by

3') if $\alpha\in R$ and $(\alpha ,\alpha )=0$ then for any
$\beta\in R$ such that $(\alpha ,\beta )\ne 0$ at least one of the
vectors $\beta + \alpha$ or $\beta - \alpha$ belongs to $R$.

This axiomatics includes the root systems of type $A(1,1)$ as well as the root systems of type $C(m,n)$ and $BC (m,n).$ We have used it in \cite{SV} to introduce a class of the deformed Calogero-Moser operators.

\section{Ring $J(\mathfrak {g})$ and supercharacters of $\mathfrak {g}$}

Let $V$ be a finite dimensional module over a basic classical Lie superalgebra $\mathfrak {g}$
with Cartan subalgebra $\mathfrak {h}.$ Let us assume for the moment that $V$ is a {\it semisimple} $\mathfrak {h}$-module, which means that $V$ can be decomposed as a sum of the one-dimensional
$\mathfrak {h}$-modules:
$$V = \bigoplus_{\lambda \in P(V)} V_{\lambda},$$
where $P(V)$ is the set of the corresponding weights $\lambda \in \mathfrak {h}^*.$
The {\it supercharacter} of $V$ is defined as
$$sch \,V = \sum_{\lambda \in P(V)} ( {\it sdim}\, V_{\lambda}) e^{\lambda},$$
where ${\it sdim}$ is the {\it superdimension} defined for any $\mathbb{Z}_{2}$-graded vector space $W= W_0\oplus W_1$ as the difference of usual dimensions of graded components:
$$ {\it sdim} \,W = {\it dim} \,W_0 - {\it dim} \, W_1.$$
By definition the supercharacter
$sch \,V \in \mathbb{Z}[\mathfrak {h}^*]$
is an element of the integral group ring of $\mathfrak {h}^*$ (considered as an abelian group).

The following proposition shows that in the context of Grothendieck ring we can restrict ourselves by the semisimple $\mathfrak {h}$-modules. First of all note that the Grothendick group has a natural basis consisting of irreducible modules. Indeed any finite dimensional module has Jordan-H\"older series, so in Grothendieck group it is equivalent to the sum of irreducible modules.

\begin{proposition}\label{semi}
Let $V$ be a finite dimensional irreducible $\mathfrak g$-module. Then $V$ is semisimple as  $\mathfrak h$-module.
\end{proposition}

\begin{proof}
Let $W\subset V$ be the maximal semisimple $\mathfrak h$-submodule. Since $V$ is finite dimensional $W$ is nontrivial.
Let us prove that $W$ is $\mathfrak g$-module. We have
$$
\mathfrak g=\mathfrak h \oplus(\bigoplus_{\alpha\in R}\mathfrak g_{\alpha}).
$$
Since $W$ is semisimple it is a direct sum of one-dimensional $\mathfrak h$-modules.
Let $w\in W$ be a generator of one of them, so that $hw=l(h)w$ for any $h \in \mathfrak h.$
Note that $xw$ is an eigenvector for $\mathfrak h$ for any $x \in \mathfrak g_{\alpha}$ since for any $h \in \mathfrak h$
$$
hxw=[h,x]w+xhw=(\alpha(h)+l(h))xw.
$$
 Now the fact that $xw$ belongs to $W$ follows from the maximality of $W$. Since $V$ is irreducible $W$ must coincide with $V.$
\end{proof}

The following general result (essentially contained in Kac \cite{Kac, Kac2}) shows that for basic classical Lie superalgebras an irreducible module is uniquely determined by its supercharacter.

\begin{proposition}
\label{char}
 Let $V,U$ be finite dimensional irreducible $\mathfrak g$-modules. If $sch \,V=sch\, U$ then  $V,U$ are isomorphic as $\mathfrak g$-modules.
\end{proposition}

\begin{proof}
By the previous proposition the modules are semisimple.
According to Kac \cite{Kac2} (see proposition 2.2) every irreducible finite dimensional module is uniquely determined by its highest weight. Since $sch\, V=sch\, U$   the modules
 $V$ and $U$ have the same highest weights and thus are isomorphic as $\mathfrak g$ modules.
\end{proof}

Now we are going to explain why the definition (\ref{Jdefin}) of the ring $J(\mathfrak{g})$ is natural in this context.

Recall that in the classical case of semisimple Lie algebras  the representation theory of $\mathfrak{sl}(2)$ plays the key role (see e.g. \cite{FH}).  In the super case it is natural to consider the Lie superalgebra $\mathfrak{sl}(1,1)$, which has three generators $H, X,Y$ ($H$ generates the even part, $X,Y$ are odd) with the following relations:
\begin{equation}
\label{sl1,1}
[H, X] = [H,Y] = [ Y,Y] = [ X,X]=0,\, [ X,Y]=H.
\end{equation}
However because of the absence of complete reducibility in the super case this Lie superalgebra alone is not enough to get the full information. We need to consider the following extension of $\mathfrak{sl}(1,1)$.
As before we will use the notation $(a)$ for the principal ideal of the integral group ring $\mathbb{Z}[\mathfrak {h}^*]$ generated by an element $a \in \mathbb{Z}[\mathfrak {h}^*].$

\begin{proposition}\label{J}
 Let  $\mathfrak{g}(\mathfrak{h},\alpha)$
be the solvable Lie superalgebra such that
$\mathfrak{g}_{0}=\mathfrak{h}$ is a commutative finite
dimensional Lie algebra, $\mathfrak{g}_{1}=Span(X, Y)$ and
the following relations hold
\begin{equation}
\label{XY}
[h,X]=\alpha(h)X,\quad [h,Y]=-\alpha(h)Y,\quad
[ Y,Y]=[ X,X]=0,\:
[ X,Y]=H,
\end{equation}
where $H \in \mathfrak{h}$ and $\alpha\ne0$ is a linear form on $\mathfrak{h}$ such that $ \alpha(H)=0.$
Then the Grothendieck ring of $\mathfrak{g}(\mathfrak{h},\alpha)$ is isomorphic to
\begin{equation}
\label{Jprop}
J(\mathfrak{g}(\mathfrak{h},\alpha))=\{ f=\sum c_{\lambda}
e^{\lambda}\mid \lambda\in\mathfrak{h}^*,\quad D_{H}f\in (e^{\alpha}-1)\},
\end{equation}
where by definition $D_{H}e^{\lambda}=\lambda(H)e^{\lambda}.$ The isomorphism is given by the supercharacter map
$
Sch: [V] \longrightarrow sch \,V.
$
\end {proposition}

\begin{proof}
Every irreducible $\mathfrak{g}(\mathfrak{h},\alpha)$-module $V$  has a unique (up to a multiple) vector $v$ such that $X v = 0, \, h v = \lambda(h)v$ for some linear form $\lambda$ on $\mathfrak{h}.$ This establishes a bijection between the  irreducible $\mathfrak{g}(\mathfrak{h},\alpha)$-modules and the elements of $\mathfrak{h}^*.$

There are two types of such modules, depending on whether $\lambda(H)=0$ or not.  In the first case the module $V=V({\lambda})$ is one-dimensional and its
supercharacter is $e^{\lambda}$. If
$\lambda(H)\ne0$ then the corresponding module $V({\lambda})$ is two-dimensional with the  supercharacter
$sch(V) = e^{\lambda}-e^{\lambda-\alpha}.$
Thus we have proved that the image of supercharacter map $Sch\, (K(\mathfrak{g}(\mathfrak{h},\alpha)))$ is contained in $J(\mathfrak{g}(\mathfrak{h},\alpha))$.

Conversely, let $f=\sum
c_{\lambda}e^{\lambda}$ belong to $J(\mathfrak{g}(\mathfrak{h},\alpha))$. By subtracting a suitable linear combination of supercharacters of the one-dimensional modules $V(\lambda)$ we can assume
that $\lambda(H)\ne0$ for all $\lambda$ from $f$. Then the condition $D_{H}f \in
(e^{\alpha}-1)$ means that
\begin{equation}
\label{H}
\sum \lambda(H) c_{\lambda} e^{\lambda}= \sum d_{\mu} (e^{\mu}-e^{\mu-\alpha}).
\end{equation}
For any $\lambda \in \mathfrak{h}^*$ define the linear functional $F_{\lambda}$ on $\mathbb{Z}[\mathfrak {h}^*]$ by
$$
F_{\lambda}(f)=\sum_{k\in \mathbb{Z}}c_{\lambda+k\alpha}.
$$
It is easy to see that the conditions $F_{\mu}(f) = 0$ for all $\mu \in \mathfrak{h}^*$
characterise the ideal $(e^{\alpha}-1).$
Applying $F_{\mu}$ to both sides of the relation (\ref{H}) and using the fact that
$\alpha(H) = 0, \, \lambda(H) \neq 0$ we deduce that $f$ itself belongs to the ideal. This means that
$f = \sum p_{\nu} (e^{\nu}-e^{\nu-\alpha})$ for some integers $p_{\nu},$ which is a linear combination of the supercharacters of the irreducible modules $V(\nu).$
\end{proof}

Any basic classical Lie superalgebra has a subalgebra isomorphic to (\ref{XY}) corresponding to any isotropic root $\alpha.$  By restricting the modules to this subalgebra we have the following

\begin{proposition}\label{image}
For any basic classical Lie superalgebra $\mathfrak{g}$ the supercharacter map $Sch$ is injective and its image $Sch\, (K(\mathfrak{g}))$ is
contained in $J(\mathfrak{g})$.
\end{proposition}

The first claim is immediate consequence of proposition \ref{char}.
The invariance with respect to the Weyl group $W_0$ follows from the fact that any finite dimensional $\mathfrak{g}$-module is also finite dimensional $\mathfrak{g}_0$-module.

This gives the proof of an easy part of the Theorem. The rest of the proof (surjectivity of the supercharacter map) is much more involved.
%The next section contains some technical results we will need.

 \section{Geometry of the highest weight set}

In this section, which is quite technical, we give the description of the set of highest weights of finite
dimensional $\mathfrak{g}$-modules in terms of the corresponding generalized root systems.
Essentially one can think of this as a geometric interpretation of the Kac conditions \cite{Kac, Kac2}.

Following Kac \cite{Kac2} we split  all basic classical Lie superalgebras $\mathfrak{g} = \mathfrak{g}_0 \oplus \mathfrak{g}_1$ into two {\it  types}, depending on
whether $\mathfrak{g}_0$-module $\mathfrak{g}_1$ is reducible (type I) or not (type II).
The Lie superalgebras $A(m,n), \, C(n)$ have type I, type II list consists of $$B(m,n), \, D(m,n),\, F(4), G(3), D(2,1, \alpha).$$
In terms of the corresponding root systems type II is characterised by the property that the even roots generate the whole dual space to Cartan subalgebra $\mathfrak{h}.$
In many respects Lie superalgebras of type II have more in common with the usual case of simple Lie algebras
than Lie superalgebras of type I. In particular, we will see that the corresponding Grothendieck rings in type II can be naturally realised as subalgebras of the polynomial algebras, while in type I it is not the case.

Let us choose a {\it distinguished  system} $B$ of simple roots in $R$, which contains only one  isotropic root $\gamma$; this is possible for any basic classical Lie superalgebra except $B(0,n),$ which has no isotropic roots
(see \cite{Kac}).
If we take away $\gamma$ from $B$ the remaining set will give the system of simple roots of the even part $\mathfrak{g}_0$ if and only if $\mathfrak{g} $ has type I.
For type II one can replace $\gamma$ in $B$ by a unique positive even root $\beta$ (called {\it special}) to get a basis of simple roots of  $\mathfrak{g}_0.$

In the rest of this section we restrict ourselves with the Lie superalgebras of type II.
The following fact, which will play an important role in our proof, can be checked case by case (see explicit formulas in the last section).

\begin{proposition}\label{fact}
For any basic classical Lie superagebra of type II except $B(0,n)$ there exists a unique decomposition of Lie algebra $\mathfrak{g}_0 = \mathfrak{g}_0^{(1)} \oplus \mathfrak{g}_0^{(2)}$ such that
the isotropic simple root $\gamma$ from distinguished  system $B$ is the difference
\begin{equation}
\label{gamma}
 \gamma = \delta - \omega
\end{equation}
of two weights $\delta$ and $\omega$ of $\mathfrak{g}_0^{(1)}$ and $ \mathfrak{g}_0^{(2)}$ respectively with the following properties:

1) $\mathfrak{g}_0^{(2)}$ is a semisimple Lie algebra and $\omega$ is its fundamental weight

2) the special root $\beta$ is a root of $\mathfrak{g}_0^{(1)}$ and $\delta = \frac{1}{2} \beta.$

\noindent In the exceptional case $B(0,n)$ we define $\mathfrak{g}_0^{(1)}= \mathfrak{g}_0$
and $\omega =0.$
\end{proposition}

\noindent {\bf Remark.}  The fundamental weight $\omega$ has the following property, which will be very important for us: it has a small orbit in the sense of Serganova (see below).

Let $\mathfrak{a}$ be a semisimple Lie algebra, $W$ be its Weyl group, which acts on the corresponding root system $R$ and weight lattice $P$ (see e.g. \cite{Bou}).
Following Serganova \cite{Serga} we call the orbit $W \omega$ of weight $\omega$ {\it small} if
for any  $x, y \in W\omega$ such that $x\ne\pm y$ the difference $x-y$ belongs to the root system $R$
of $\mathfrak{a}$. Such orbits play a special role in the classification of the generalized root systems.

Let $\mathfrak{g}$ be a basic classical Lie superalgebra of type II,
$\mathfrak{a}= \mathfrak{g}_0^{(2)}$ and $ \omega$ as in Proposition \ref{fact}.
Define a positive integer $k = k (\mathfrak{g})$ as
\begin{equation}
\label{k}
k=
\frac{1}{2} \mid W\omega\mid, 
\end{equation}
where $W$ is the Weyl group of $\mathfrak{a}$ and $\mid W\omega\mid$ is the number of elements in the orbit of the weight $\omega.$

For any positive integer $j\le k$ consider a subset $L_j \subset P$ of the weight lattice of $\mathfrak{a}$ defined by the relations
\begin{equation}
\label{Lj}
 F(\nu)\ne0,\: F(\nu-\omega)=0,\dots,
F(\nu-(j-1)\omega)=0,\,
(\nu,\omega)=(\rho+(j-k)\omega,\omega),
\end{equation}
where $$F(\nu)=\prod_{\alpha\in R^+}(\nu,\alpha)$$ and $\rho$ is the half of the sum of positive roots $\alpha\in R^+$ of $\mathfrak{a}.$
In particular,
$$L_1 = \{\nu \in P \mid  F(\nu)\ne0, (\nu,\omega)=(\rho+(1-k)\omega,\omega)\}.$$

Let $\Lambda$ be a highest weight of Lie algebra $\mathfrak{g}_0$
and $\lambda$ be its projection  to the
weight lattice of $\mathfrak{a}= \mathfrak{g}_0^{(2)}$ with respect to the decomposition $\mathfrak{g}_0 = \mathfrak{g}_0^{(1)} \oplus \mathfrak{g}_0^{(2)}.$
%\begin{equation}
%\label{rho}
%\rho
%\end{equation}
Define an integer $ j(\Lambda)$ by the formula
\begin{equation}
\label{j}
j({\Lambda})=k-\frac{(\Lambda,\delta)}{(\delta,\delta)}
\end{equation}
where $\delta$ is the same as in Proposition \ref{fact}.
This number was implicitly used by Kac in \cite{Kac2}.

Define the following set  $X (\mathfrak{g})$ consisting of
the highest weights $\Lambda$ of $\mathfrak{g}_0$  such that either $j(\Lambda) \leq 0$ or the $W$-orbit of $\lambda + \rho$ intersects the set $L_j$ for some $j=1, \dots,k.$  

The main result of this section is the following 

\begin{thm}\label{kaccond}
For any basic classical Lie superalgebra $\mathfrak{g}$ of type II the set $X (\mathfrak{g})$ coincides with the set of the highest weights of the finite dimensional representations of $\mathfrak{g}.$
\end{thm}

The rest of the section is the proof of this theorem.
Let us define the {\it support} $Supp(\varphi)$ of an element $\varphi \in \mathbb{Z}[ P]$ as the set of weights $\nu \in P$ in the representation $\varphi = \sum \varphi(\nu) e^{\nu},$ for which $\varphi(\nu)$ is not zero. Define also the {\it alternation} operation on $\mathbb{Z}[ P]$ as
\begin{equation}
\label{Alt}
Alt(\varphi)=\sum_{w\in W}\varepsilon(w){w(\varphi)},
\end{equation}
where by definition $w(e^{\nu}) = e^{w\nu}$ and $\varepsilon: W \rightarrow \pm1$ is the sign homomorphism.
\begin{lemma}\label{alt}
 Let
$\omega$ be a weight such that the orbit $W\omega$ is small.
Consider
$\varphi \in \mathbb{Z}[ P]$
such that $Alt(\varphi)=0,$
$Supp(\varphi)$ is contained in the hyperplane $(\nu,\omega)=c$ for some $c$
and for every $\nu \in Supp(\varphi)$
$$F(\nu)=\prod_{\alpha\in R^+}(\nu,\alpha)\neq 0.$$
Then

1) if $c\ne0$ then for any $t \in \mathbb{Z}$ \, $Alt(\varphi
e^{t\omega})=0;$

2) if $c=0$ then the same is true if there exists $\sigma_{0}\in W$ such that $\sigma_{0}\omega=-\omega$ and
$\sigma_{0}\varphi=\varphi,\: \varepsilon(\sigma_{0})=1.$
\end{lemma}

\begin{proof}
We have
 $$
 Alt(\varphi)=\sum_{(\nu,\omega)=c}\varphi(\nu)Alt(e^{\nu})=0.
 $$
 Since $F(\nu) \neq 0$ the elements $Alt(e^{\nu})$ are non-zero and linearly independent for $\nu$ from different orbits of $W.$ Thus the last equality is equivalent to
\begin{equation}
\label{sigma}
\sum_{\sigma\in W}\varepsilon(\sigma)\varphi(\sigma \nu)=0
\end{equation}
for any $\nu$ from the support of $\varphi$.

Let $c\ne0.$ Fix $\nu \in Supp(\varphi)$ and consider $\sigma \in W$ such that  $\varphi(\sigma \nu) \neq 0,$ in particular $(\sigma\nu,\omega)=c$. We have $(\nu,\omega-\sigma^{-1}\omega)=0,\:(\nu,\omega+\sigma^{-1}\omega)=2c\ne0$. Since the orbit of $\omega$ is small and $F(\nu)\neq 0$ this implies that
 $\omega=\sigma^{-1}\omega$ and  therefore $\sigma$ belongs to the stabiliser $W_{\omega} \subset W$ of $\omega$.
 Thus the relation (\ref{sigma}) is equivalent to
% $$
% \sum_{\sigma\in W_{\omega}}\varepsilon(\sigma)\varphi(\sigma \nu)=0,
 %$$
%which can be rewritten as
 $$
 \sum_{\sigma\in W_{\omega}}\varepsilon(\sigma)\sigma(\varphi)=0.
 $$
Since  $\omega$ is invariant under $W_{\omega}$ this implies
$$
\sum_{\sigma\in W_{\omega}}\varepsilon(\sigma){\sigma(\varphi e^{t\omega})}=0
$$
and thus
$$
\sum_{\sigma\in W}\varepsilon(\sigma){\sigma(\varphi e^{t\omega})}=0.
$$
This proves the first part.

When $c=0$ similar arguments lead to the relation
 $$
 \sum_{\sigma\in W_{\pm\omega}}\varepsilon(\sigma)\sigma(\varphi)=0,
 $$
where $W_{\pm\omega}$  is the stabiliser of the set $\pm\omega.$  From the conditions of the lemma it follows that  $W_{\pm\omega}$ is generated by $W_{\omega}$ and $\sigma_{0}.$
Since $\varepsilon (\sigma_{0})=1$ and $ \sigma_{0}\varphi=\varphi$ we can replace in this last formula  $W_{\pm\omega}$ by $W_{\omega}$ and repeat the previous arguments to complete the proof. 
\end{proof}

%Let $\mathfrak{g}$ be one of the following semisimple Lie algebras from the list
%\begin{equation}
%\label{list}
%A(n),\:B(n),\:C(n),\:D(n),\:G(2),\:B(3),\:A(1)\oplus A(1)
%\end{equation}
%and $\omega$ be the
%following fundamental weights of the corresponding standard representations
%\begin{equation}
%\label{list1}
%\varepsilon_{1},\:\varepsilon_{1},\:\varepsilon_{1},\:\varepsilon_{1},\:\varepsilon_{2}+\varepsilon_{3},
%\frac{1}{2}(\varepsilon_{1}+\varepsilon_{2}+\varepsilon_{3}),\:\varepsilon_{2}+\varepsilon_{3}
%\end{equation}
%respectively.

%The following lemma will play a key role.

Recall that for any $\omega \in P$ the derivative $D_{\omega}$ is determined by the relation  $D_{\omega} e^{\lambda} = (\omega,\lambda) e^{\lambda}.$ The condition that $D_{\omega} \varphi = 0$ is equivalent to the support of $\varphi$ to be
contained in the hyperplane $(\omega,\lambda)=0.$

\begin{lemma}\label{small}
Let $\mathfrak{g}$ be a basic classical Lie superalgebra,
$\mathfrak{a} = \mathfrak{g}_0^{(2)}$ and $ \omega$ as in Proposition \ref{fact}, $k$ defined by (\ref{k}), $W$ be the Weyl group of $\mathfrak{a}$ acting on the corresponding weight lattice $P$.
Consider a function of the form
\begin{equation}
\label{smal}
\varphi=
\sum_{i=1}^{k}(e^{(k-i)\omega}+e^{-(k-i)\omega})f_{i},
\end{equation}
where $f_{i} \in \mathbb{Z}[ P]^{W} $ are some exponential $W$-invariants. Suppose that
$D_{\omega} \varphi = 0$ and consider the first non-zero coefficient $f_{j}$ in $\varphi$
(so that $f_{1}=f_{2}=\dots=f_{j-1}=0$ for some $j\le k$).

Then $f_{j}$ is a linear
combination of the characters of irreducible representations  of $\mathfrak a$ with
the highest weights $\lambda$ such that  the orbit $W(\lambda+\rho)$ intersects the set
$L_j$ defined above.
\end{lemma}

\begin{proof}
Since $D_{\omega} \varphi = 0$ the support of $\varphi$ is contained in the hyperplane
$(\omega,\mu)=0$. We can write the function $\varphi$ as the sum
$\varphi= \varphi_{1}+\dots+\varphi_{j}+\psi_{j},$ where the support of
$\varphi_{j}e^{\rho+(j-k)\omega}$ is contained in $L_j$ and the support of $\psi_{j}e^{\rho+(i-k)\omega} $ is not contained
in $L_i$ for all $i=1,\dots, j$. Let us multiply (\ref{smal}) consequently by $e^{\rho+(1-k)\omega},\:e^{\rho+(2-k)\omega},\:\dots,
\:e^{\rho+(j-k)\omega}\:$ and then apply the alternation operation (\ref{Alt}). Then
from the definition of the sets $L_j$ we have
$$
Alt(e^{\rho})f_{1}=Alt(\varphi_{1}e^{\rho+(1-k)\omega}),
$$
$$
Alt(e^{\rho+\omega})f_{1}+Alt(e^{\rho})f_{2}=Alt(\varphi_{1}e^{\rho+(2-k)\omega})+
Alt(\varphi_{2}e^{\rho+(2-k)\omega}),
$$
$$
\dots\dots\dots
$$
$$
Alt(e^{\rho+(j-1)\omega})f_{1}+\dots+Alt(e^{\rho})f_{j}=Alt(\varphi_{1}e^{\rho+(j-k)\omega})+
\dots+ Alt(\varphi_{j}e^{\rho+(j-k)\omega}).
$$
Suppose that $f_{1}=f_{2}=\dots=f_{j-1}=0$. Then from the
first equation we see that
$Alt(\varphi_{1}e^{\rho+(1-k)\omega})=0$. One can verify that $(\rho+(j-k)\omega, \omega)=0$ if and only if $j=1$ and $(\mathfrak a,\omega)$ must be either $(D(m),\varepsilon_{1})$ or $(B_{3},1/2(\varepsilon_{1}+\varepsilon_{2}+\varepsilon_{3}).$  In both of these cases  we can find $\sigma_{0}$ such that $\varepsilon(\sigma_{0})=1,\:\sigma_{0}(\omega)=-\omega,\: \sigma_{0}\varphi=\varphi$, so we can apply Lemma \ref{alt} to show that
$Alt(\varphi_{1}e^{\rho-(k-i)\omega})=0$ for $i=1,\dots,j$. Similarly
from the second equation
$Alt(\varphi_{2}e^{\rho+(2-k)\omega})=0$ and by applying again Lemma \ref{alt} we have
$Alt(\varphi_{2}e^{\rho+(i-k)\omega})=0$ for $i=2,\dots,j$ and eventually
$$
Alt(e^{\rho})f_{j}=Alt(\varphi_{j}e^{\rho+(j-k)\omega}).
$$
Now the claim follows from the classical {\it Weyl character formula}  (see e.g. \cite{Serre}) for the representation with highest weight $\lambda$:
\begin{equation}
\label{Weyl}
ch \, V^{\lambda} = \frac {Alt (e^{\lambda + \rho})}{Alt (e^{\rho})}.
\end{equation}
\end{proof}

Now we need the conditions on the highest weights of the finite dimensional representations, which were found by Kac \cite{Kac}.
In the following Lemma, which is a reformulation of proposition 2.3 from \cite{Kac}, we use the basis of the weight lattice of $\mathfrak{g}_0$ described in Section 7.

\begin{lemma}\label{Kac} (Kac \cite{Kac}).
For the basic classical Lie superalgebras $\mathfrak{g}$ of type II a highest weight $\nu$ of $\mathfrak{g}_0$ is a highest weight of finite dimensional irreducible $\mathfrak{g}$-module if and only if one  of the corresponding conditions is satisfied:

1) $\mathfrak{g}=B(m,n)$, \, $\Lambda =  (\mu_1, \dots, \mu_n, \lambda_1,\dots, \lambda_m)$
\begin{itemize}
\item $\mu_{n}\ge m$

\item $\mu_{n}=m-j,\: 0<j\le m $ and
$\lambda_{m}=\lambda_{m-1}=\dots=\lambda_{m-j+1}=0$
\end {itemize}

2) $\mathfrak{g}=D(m,n)$, \, $\Lambda = (\mu_1, \dots, \mu_n, \lambda_1,\dots, \lambda_m)$

\begin{itemize}
\item $\mu_{n}\ge m$

\item $\mu_{n}=m-j,\: 0<j\le m $ and
$\lambda_{m}=\lambda_{m-1}=\dots=\lambda_{m-j+1}=0$
\end {itemize}

 3) $\mathfrak{g}=G(3)$, \,$ \Lambda = (\mu, \lambda_1, \lambda_2) $

\begin{itemize}
\item $\mu\ge 3$

\item $\mu=2,\:\lambda_{2}=2\lambda_{1}$

\item $\mu=0,\, \lambda_{1}=\lambda_{2}=0$

\end {itemize}

4) $\mathfrak{g}=F(4)$, \, $\Lambda = (\mu, \lambda_1, \lambda_2, \lambda_3)$

\begin{itemize}
\item $\mu\ge 4$

\item $\mu=3,\:\lambda_{1}=\lambda_{2}+\lambda_{3}-1/2$

\item $\mu=2,\:\lambda_{1}=\lambda_{2},\:\lambda_{3}=0$

\item $\mu=0,\:\lambda_{1}=\lambda_{2}=\lambda_{3}=0$
\end {itemize}

5) $\mathfrak{g}=D(2,1,\alpha)$, \, $\Lambda = (\lambda_1, \lambda_2, \lambda_3)$

\begin{itemize}

\item $\lambda_{1}\ge 2$

\item  $\lambda_{1}=1$ , $\alpha$ is rational and
$\lambda_{2}-1=|\alpha|(\lambda_{3}-1)$

\item $\lambda_{1}=0,\lambda_{2}=\lambda_{3}=0.$
\end {itemize}
\end{lemma}

%For the Lie superalgebras of the type I  we will derive our main Theorem from the %explicit description of the ring $J(\mathfrak{g})$ given in Section 7.

Now we are ready to prove Theorem 5.2.

Let $\mathfrak{g}$ be a basic classical Lie superalgebra of type II,
$\mathfrak{a}, \, \omega, \, k, \, W$ be the same as in Lemma \ref{small},
$\Lambda$ be a highest weight of Lie algebra $\mathfrak{g}_0$, $\lambda$ be its projection to the
weight lattice of $\mathfrak{a}$ and $j = j(\Lambda)$ is defined by the formula
(\ref{j}).

We are going to show that the conditions defining the set $X (\mathfrak{g})$ are equivalent to the Kac's conditions from the previous Lemma. First of all an easy check shows that  in each case the condition $j(\Lambda) \leq 0$ is equivalent to the first of Kac's conditions. 

Let us consider now the condition that $W(\lambda + \rho)$ intersects the set $L_j.$
We will see that in that case $j = j(\Lambda).$

By definition $L_j$ is described by the following system for the weights $\nu$ of $\mathfrak{a}$
%$$
%\left\{ \begin{array}{r}
%(\mu,\omega)=0,\\
%F(\mu+\rho-(k-1)\omega)=0\\
%\dots\\
%F(\mu+\rho-(k-j+1)\omega)=0\\
%F(\mu+\rho-(k-j)\omega)\ne0\\
%\end{array}
%\right.
%$$
%Let $\nu=\mu+\rho-(k-j)\omega$, then the previous system is
%equivalent to the following
$$
\left\{ \begin{array}{r}
F(\nu)\ne0\\
F(\nu-\omega)=0\\
F(\nu-2\omega)=0\\
\dots\\
F(\nu-(j-1)\omega)=0\\
(\rho-(k-j)\omega,\omega)=(\nu,\omega).\\
\end{array}
\right.
$$
Consider this system in each case separately.

 1) Let $\mathfrak{g}=B(m,n)$ with $m>0$, $\mathfrak{a}=B(m)$, then $k=m, \, \omega=\varepsilon_{1}, \,\rho = \sum_{i=1}^m (m-i+1/2) \varepsilon_i$ and
 $$
 F(\nu)=\prod_{p=1}^m\nu_{p}\prod_{p<q}(\nu_{p}^2-\nu_{q}^2).
 $$
Since $F(\nu) \neq 0$ all $\nu_i$ are non-zero and pairwise different.
 The condition that $(\rho-(k-j)\omega,\omega)=(\nu,\omega)$ means that $\nu_{1}=j-1/2$.
 Then we have the following system
 $$
 \left\{ \begin{array}{r}
(\nu_{2}^2-(j-3/2)^2)(\nu_{3}^2-(j-3/2)^2)\dots(\nu_{m}^2-(j-3/2)^2)=0\\
(\nu_{2}^2-(j-5/2)^2)(\nu_{3}^2-(j-5/2)^2)\dots(\nu_{m}^2-(j-5/2)^2)=0\\
\dots\\
(\nu_{2}^2-(1/2)^2)(\nu_{3}^2-(1/2)^2)\dots(\nu_{m}^2-(1/2)^2)=0.\\
\end{array}
\right.
 $$
 The first equation implies that one of $\nu_i$ equals to $j-3/2$, the second one implies that one of them is $j-5/2$ and so on.
So if $W(\lambda+\rho)\cap L_j\ne\emptyset $ then
$\lambda_{m}=\lambda_{m-1}=\dots=\lambda_{m-j+1}=0$, which is one of the corresponding conditions in Lemma \ref{Kac}.
In the case $B(0,n)$ we have the only condition $j(\Lambda) \leq 0,$ which is equivalent to $\mu_n \geq 0.$

2) When $\mathfrak{g}=D(m,n),\, \mathfrak{a}=D(m)$ we have $k =m,\, \omega=\varepsilon_{1}, \, \rho = \sum_{i=1}^m (m-i) \varepsilon_i$ and
 $$
 F(\nu)=\prod_{p<q}(\nu_{p}^2-\nu_{q}^2).
 $$
In that case  the condition $(\rho-(k-j+1)\omega,\omega)=(\nu,\omega)$ implies that $\nu_{1}=j-1$ and we have the following system
 $$
 \left\{ \begin{array}{r}
(\nu_{2}^2-(j-2)^2)(\nu_{3}^2-(j-2)^2)\dots(\nu_{m}^2-(j-2)^2)=0\\
(\nu_{2}^2-(j-3)^2)(\nu_{3}^2-(j-3)^2)\dots(\nu_{m}^2-(j-3)^2)=0\\
\dots\\
\nu_{2}^2\nu_{3}^2\dots\nu_{m}^2=0.\\
\end{array}
\right.
 $$
If $W(\lambda+\rho)\cap L_j\ne\emptyset $ we have similarly again
$\lambda_{m}=\lambda_{m-1}=\dots=\lambda_{m-j+1}=0$.

3) Let $\mathfrak{g}=G(3), \, \mathfrak{a}=G(2)$, then
$k=3,\,\omega=\varepsilon_{1}+\varepsilon_{2},\, \rho = 2 \varepsilon_1 + 3 \varepsilon_2$ and
 $$
 F(\nu)=\nu_{1}\nu_{2}(\nu_{2}-\nu_{1})(\nu_{1}+\nu_{2})(2\nu_{1}-\nu_{2})(2\nu_{2}-\nu_{1}).
 $$
 The condition $(\rho-(k-j)\omega,\omega)=(\nu,\omega)$ means that
 $\nu_{1}+\nu_{2}=2j-1$. If $j=1$ we have
$
\nu_{1}+\nu_{2}=1,\, F(\nu)\ne0.
$
One check that in that case $\nu\in W(\lambda+\rho)$ only if $\lambda_{2}=2\lambda_{1}.$

If $j=2$ we have the conditions
$
\nu_{1}+\nu_{2}=3,\, 
F(\nu_{1}-1,\nu_{2}-1)=0,\ F(\nu)\ne0,
 $
which can not be satisfied for $\nu \in W(\lambda+\rho)$.

If $j=3$ we have
$
\nu_{1}+\nu_{2}=3,\,
F(\nu_{1}-1,\nu_{2}-1)=0,\,
F(\nu_{1}-2,\nu_{2}-2)=0,\,
F(\nu)\ne0,
$
which imply that if $\nu\in W(\lambda+\rho)$ then $\lambda_{2}=\lambda_{1}=0$
in agreement with Lemma \ref{Kac}.

 4) If $\mathfrak{g}=F(4),\, \mathfrak{a}=B(3)$ then
$k=4,\,\omega=\frac{1}{2}(\varepsilon_{1}+\varepsilon_{2}+\varepsilon_{3}), \rho = \frac{5}{2} \varepsilon_1+\frac{3}{2}  \varepsilon_2 +\frac{1}{2}  \varepsilon_3,$
 $$
 F(\nu)=\nu_{1}\nu_{2}\nu_{3}(\nu^2_{1}-\nu^2_{2})
 (\nu^2_{1}-\nu^2_{3})(\nu^2_{2}-\nu^2_{3}).
 $$
 The condition $(\rho-(3-j)\omega,\omega)=(\nu,\omega)$ means that
 $\nu_{1}+\nu_{2}+\nu_{3}=\frac{3}{2}(j-1)$. If $j=1$ we have the conditions
$\nu_{1}+\nu_{2}+\nu_{3}=0,\,
F(\nu)\ne0,
 $
 which imply that if $\nu\in W(\lambda+\rho)$ then
$\lambda_{1}=\lambda_{2}+\lambda_{3}-1/2.$

If $j=2$ we have
$
\nu_{1}+\nu_{2}+\nu_{3}=3/2,\,
F(\nu_{1}-1/2,\nu_{2}-1/2,\nu_{3}-1/2)=0,\,
F(\nu)\ne0.
$
One can check that  if $\nu\in W(\lambda+\rho)$ then
$\lambda_{1}=\lambda_{2},\:\lambda_{3}=0.$

If $j=3$ we have the conditions
 $
\nu_{1}+\nu_{2}+\nu_{3}=3,\,
F(\nu_{1}-1/2,\nu_{2}-1/2,\nu_{3}-1/2)=0,\,
F(\nu_{1}-1,\nu_{2}-1,\nu_{3}-1)=0,\,
F(\nu)\ne0,
 $ which can not be satisfied for $\nu \in W(\lambda+\rho)$.

If $j=4$ we have
$
\nu_{1}+\nu_{2}+\nu_{3}=9/2,\,
F(\nu_{1}-1/2,\nu_{2}-1/2,\nu_{3}-1/2)=0,\,
F(\nu_{1}-1,\nu_{2}-1,\nu_{3}-1)=0,\,
F(\nu_{1}-3/2,\nu_{2}-3/2,\nu_{3}-3/2)=0,\,
F(\nu)\ne0.
$
In that case $\nu\in W(\lambda+\rho)$ only if
$\lambda_{1}=\lambda_{2}=\lambda_{3}=0.$

5) Let $\mathfrak{g}=D(2,1,\alpha),\, \mathfrak{a}= A_1 \oplus A_1$, then
$k=2,\, \omega=\varepsilon_{2}+\varepsilon_{3}, \rho = -\varepsilon_{2}-\varepsilon_{3}$ and
 $
 F(\nu)=4\alpha \nu_{1}\nu_{2}.
 $
 The condition $(\rho-(2-j)\omega,\omega)=(\nu,\omega)$ means that
 $\nu_{1}+\alpha\nu_{2}=j-1$. If $j=1$ we have
$
\nu_{1}+\alpha\nu_{2}=0,\,
\nu_{1}\nu_{2}\ne0.
$
If $\alpha$ is irrational the system has no integer solution. If
$\alpha$ is rational and $\nu\in W(\lambda+\rho)$ then
$(\lambda_{1}+1)=|\alpha|(\lambda_{2}+1).$

If $j=2$ the conditions
$
\nu_{1}+\nu_{2}=1,\,
(\nu_{1}-1)(\nu_{2}-1)=0,\,
\nu_{1}\nu_{2}\ne0,
$
imply that if $\nu\in W(\lambda+\rho)$ then $\lambda_{1}=\lambda_{2}=0$
in agreement with Lemma \ref{Kac}.

This completes the proof of Theorem 5.2.\footnote{As we have recently learnt from Serganova a different description of the set of highest weights can be found in \cite{Serga4}.}

\section{Proof of the main Theorem}

Let $\mathfrak g$ be  a basic classical  Lie superalgebra of type II,
 $\mathfrak{g}_0 = \mathfrak{g}_0^{(1)} \oplus \mathfrak{g}_0^{(2)}$ be the decomposition of the corresponding Lie algebra  $\mathfrak{g}_0$ from Proposition \ref{fact}, $\gamma = \delta -\omega$ be the same as in (\ref{gamma}). The root system $R_{0}$ of $\mathfrak{g}_0$ is a disjoint union $R_{0}^{(1)}\cup R_{0}^{(2)}$ of root systems  of $\mathfrak{g}_0^{(1)}$ and $\mathfrak{g}_0^{(2)}$.

Let us introduce the following {\it partial order} $\succ$ on the weight lattice $P(R_{0}^{(1)})$: we say that $\mu \succeq 0$ if and only if $\mu$ is a sum of simple roots from $R_{0}^{(1)}$ and the weight $\delta$ with nonnegative integer coefficients.

\begin{lemma}\label{schar}
 Let $V^{\Lambda}$ be an irreducible finite dimensional $\mathfrak g$-module with highest weight $\Lambda$
and $\mu,\lambda $ be the projections of $\Lambda$ on $P(R_{0}^{(1)})$ and  $P(R_{0}^{(2)})$ respectively. Then
the supercharacter of $V^{\Lambda}$ can be represented as
\begin{equation}
\label{schara}
sch(V^{\Lambda})=e^{\mu}ch(V^{\lambda})+\sum_{\tilde\mu\prec\mu}e^{\tilde\mu}F_{\tilde\mu},\quad F_{\tilde\mu}\in \mathbb{Z}[P(R_{0}^{(2)})],
\end{equation}
where $\prec$ means partial order introduced above and $V^{\lambda}$ is the irreducible $\mathfrak{g}_0^{(2)}$-module with highest weight $\lambda.$
\end{lemma}

\begin{proof}
Consider $V^{\Lambda}$ as $\mathfrak{g}_0^{(1)}$-module and introduce the subspace $W \subset V^{\Lambda}$ consisting of all vectors of weight $\mu.$
Let us prove that $W$ as a module over Lie algebra $\mathfrak{g}_0^{(2)}$ is irreducible.
It is enough to prove that it is a highest  weight module over $\mathfrak{g}_0^{(2)}$.
Let $v\in W$ be a vector of weight $\tilde\Lambda$ with respect to the Cartan subalgebra $\mathfrak h$ of $\mathfrak g$.
Then $v=uv_{\Lambda},$  where $v_{\Lambda}$ is the highest weight vector of  $V^{\Lambda}$ and
$u$ is a linear combination of the elements of the form
$$
\prod_{\alpha\in (R_{0}^{1})^+}X_{-\alpha}^{n_{\alpha}}\prod_{\gamma\in R_{1}^+}X_{-\gamma}^{n_{\gamma}}\prod_{\beta\in (R_{0}^{2})^+}X_{-\beta}^{n_{\beta}},
$$
where $R_1$ is the set of roots of $\mathfrak g_1$ and $X_{\alpha}$ is an element from the corresponding root subspace of $\mathfrak g.$
We have
$$
\Lambda-\tilde\Lambda=\sum_{\alpha\in (R_{0}^{1})^+}n_{\alpha} \alpha+\sum_{\gamma\in R_{1}^+}{n_{\gamma}} \gamma+\sum_{\beta\in (R_{0}^{2})^+}{n_{\beta}}\beta.
$$
Let $\tilde\mu$ be the projection of $\tilde\Lambda$ to $P(R_{0}^{(1)}).$ It is easy to check case by case that the condition $\tilde\mu=\mu$ implies
$n_{\alpha}=n_{\gamma}=0$ for any $\alpha\in (R_{0}^{1})^+,\:\gamma\in R_{1}^+$.
This proves the irreducibility of $W$ and justifies the first term in the right hand side of (\ref{schara}).

To prove the form of the remainder in the formula (\ref{schara}) we note that if $\tilde\mu\neq\mu$ then $\tilde\Lambda < \Lambda$ with respect to the partial order defined by $R^+$ and hence
$\tilde\mu\prec \mu$ with respect to the partial order defined above. Lemma is proved.
\end{proof}

The following key lemma establishes the link between the ring $J(\mathfrak{g})$ and
the supercharacters of $\mathfrak{g}$.

\begin{lemma}\label{highest}
Consider any
$
f=\sum_{\mu}e^{\mu}F_{\mu}\in J(\mathfrak{g}),
$
 $\mu\in P(R_{0}^{(1)}),\:F_{\mu}\in \mathbb{Z}[P(R_{0}^{(2)})]$.
 Let $\mu_*$ be a maximal  with respect to the partial order $\succ$ among all $\mu$ such that $F_{\mu}\ne0$ and
 $j = j(\mu_*)$ be defined by the formula (\ref{j}).

If $j >0$ then $F_{\mu_*}$ is a linear
combination of the characters of irreducible representations  of $\mathfrak a=\mathfrak{g}_0^{(2)}$ with
the highest weights $\lambda$ such that  the orbit $W(\lambda+\rho)$ intersects the set
$L_j$ defined by (\ref{Lj}).
\end{lemma}

\begin{proof}
Since $\mu_*$ is maximal with respect to partial order $\succ$ it is also maximal with respect to the partial order defined by $(R_{0}^{(1)})^+$. Because of the symmetry  of $f$ with respect to the Weyl group of the root system $R_{0}^{(1)}$ the weight  $\mu_*$ is dominant. From the definition of the ring $J(\mathfrak{g})$ we have
$$
D_{\gamma}\left(\sum_{\mu^{\perp}=\mu_*^{\perp}}e^{\mu}F_{\mu}\right)\in(e^{\gamma}-1),
$$
where $\gamma$ is the same as in Proposition \ref{fact} and $\mu^{\perp}$ is the component of $\mu$ perpendicular to $\delta.$
This can be rewritten as
\begin{equation}
\label{D}
D_{\gamma} \phi \in(e^{\gamma}-1),
\end{equation}
where $$\phi = \sum_{\mu^{\perp}=\mu_*^{\perp}}\left(e^{\frac{(\mu,\delta)}{(\delta,\delta)}\delta}+e^{-\frac{(\mu,\delta)}{(\delta,\delta)}\delta}\right)F_{\mu}$$ (we have used the symmetry with respect to the root $2\delta$).

Let $\varphi$ be the restriction of $\phi$ on the hyperplane $\gamma = 0,$ where we consider weights as linear functions on Cartan subalgebra $\mathfrak h.$
Using the relation $\gamma = \delta - \omega$ we can rewrite (\ref{D}) as
$
D_{\omega} \varphi = 0.
$
The conditions $\mu \prec \mu_*, \, \mu^{\perp}=\mu_*^{\perp}$ imply that $\frac{(\mu_*,\delta)}{(\delta,\delta)} > \frac{(\mu,\delta)}{(\delta,\delta)}.$  We have 
$$
\varphi=(e^{(k-j)\omega}+e^{-(k-j)\omega})
F_{\mu_*}+\sum_{0\le  l < k-j}(e^{l\omega}+e^{-l\omega})F_{l},\,\,
j = k - \frac{(\mu_*,\delta)}{(\delta,\delta)}.$$
Since $F_{\mu_*},F_{l}$ are invariant  with respect to the Weyl group of $ R_{0}^{(2)}$ for $ \:0\le l < k-j$ we can apply now Lemma \ref{small} to conclude the proof.
% The same arguments can be applied for the case $-\omega\notin W\omega$.
% The claim now follows from the Lemma \ref{Kac1}.
 \end{proof}

Now we are ready to prove our main Theorem from the Introduction for the basic simple Lie superalgebras of type II.

Consider any element $f\in J(\mathfrak{g})$ and write it as in Lemma \ref{highest}
in the form
$
f=\sum_{\mu}e^{\mu}F_{\mu},
$
where $\mu\in P(R_{0}^{(1)}),\:F_{\mu}\in \mathbb{Z}[P(R_{0}^{(2)})]$.
Let $H(f) = \{\mu_1, \dots, \mu_N\}$ be the set consisting of all the maximal elements among all $\mu$ such that $F_{\mu}\ne0$ with respect to the partial order introduced above. Let $S(f)$ be the finite set of highest weights of the Lie algebra $\frak g_{0}^{(1)}$ which are less or equal than some of $\mu_i$ from $H(f) $ and $M=M(f)$ be the number of elements in the set $S(f)$.

According to Theorem \ref{kaccond} and Lemmas \ref{schar}, \ref{highest}
there are irreducible finite dimensional $\mathfrak{g}$-modules $V^{\Lambda_{1}}, \dots, V^{\Lambda_{K}}$ and integers $n_{1},\dots,n_{K}$ such that
 $$
\tilde f = f-\sum_{l=1}^K n_{l}sch(V^{\Lambda_{l}})=\sum e^{\tilde\mu}\tilde F_{\tilde\mu},
 $$
where in the last sum 
all $\tilde\mu$ are strictly less than some of $\mu_i.$ In particular this implies that none of $\mu_i$ belongs to $S(\tilde f) \subset S(f)$ and therefore $M(\tilde f) < M(f).$
Induction in $M$ completes the proof of the Theorem for type II.

\medskip

{\bf Example.} Let us illustrate the proof in the case of $G(3)$.
In this case $\frak{g}_{0}= \frak{g}^{(1)}_{0}\oplus \frak{g}^{(2)}_{0}$, where $\frak{g}^{(1)}_{0}=\frak{sl}(2),\,\ \frak{g}^{(2)}_{0}=G(2).$ Therefore  $ P(R_{0}^{(1)})=\Bbb Z$ and  the partial order introduced above coincides  with the natural order on $\Bbb Z$.  We have $\gamma=\delta-\omega$, where $\delta$ is  the  only fundamental  weight of  $\frak{sl}(2)$ and $\omega$ is  the second fundamental weight of $G(2)$. Thus for any $f\in J(\mathfrak{g}) $ the set $H(f)$ contains only one element $l\delta$ with some integer $l \geq 0$  and the corresponding $M(f) =l+1$.

%\medskip

 To prove the Theorem for type I we use the explicit description of the ring $J(\mathfrak
{g})$ given in the next section and the following notion of Kac module.

If $\frak g$ is a basic classical Lie superalgebra of type I then  $\frak g_{0}$-module 
$\frak g_{1}$ is a direct sum of two irreducible modules $\frak g_{1}=\frak g_{1}^{+}\oplus \frak g_{1}^{-}$, where $\frak g_{1}^{-}$ is linearly  generated by negative odd roots and $\frak g_{1}^{+}$ is linearly  generated by positive odd roots. One can check that $\frak g_{0}\oplus\frak g_{1}^{+}$ is a subalgebra of $\frak g$, so for every irreducible finite-dimensional $\frak g_{0}$ module $V_{0}$ we can define {\it Kac module} 
$$
K(V_{0})=U(\frak g)\otimes_{U(\frak g_{0}\oplus\frak g_{1}^{+})}V_{0}
$$
where $\frak g_{1}^{+}$ acts trivially on $V_{0}$ (see \cite{Kac}). 
Kac module is a finite-dimensional analogue of Verma module. Namely, if $\lambda$ is the highest weight of $V_{0}$, then every finite-dimensional $\frak g$-module with the same highest weight $\lambda$ is the quotient of $K(V_{0})$. It is easy to see that the character of Kac module can be given by the following formula
\begin{equation}
\label{kacm}
schK(V_{0})=\prod_{\alpha\in R_{1}^+}(1-e^{-\alpha})chV_{0}.
\end{equation}

Let us proceed with the proof now. 

Consider first the case $A(n,m)$ with $ m \neq n$. The corresponding ring  $J(\mathfrak
{g})$ is described by Proposition \ref{amnneq} and can be represented as a sum $J(\mathfrak
{g}) = \bigoplus_{a\in\mathbb C/\mathbb Z}J(\mathfrak{g})_{a}$.
 Comparing formulae (\ref{kacm}) and (\ref{pr1}) we see that the components $J(\mathfrak{g})_{a}$ with $a \notin {\mathbb Z}$   are spanned over $\mathbb Z$ by the supercharacters of Kac modules. According to the last statement of proposition \ref{amnneq} the component $J(\mathfrak{g})_{0}$ is generated over $\mathbb Z$ by $h_k$ and $h^*_k,$ which are the supercharacters of $k$-th symmetric power of the standard representation and its dual.
This proves the theorem in this case.

In the $A(n,n)$ case with $n \neq 1$ according to Proposition \ref{ann} the ring $J(\mathfrak
{g})_{0} $ is spanned over $\mathbb Z$ by the products $h_1^{m_1} h_2^{m_2}\dots h_1^{*n_1}h_2^{*n_2}\dots$ with the condition that the total degree $m_1+ 2m_2 +\dots- n_1 - 2n_2-\dots$ is equal to $0.$ It is easy to see that if $V$ is the standard representation of $\mathfrak{gl}(n+1,n+1)$ such a product is the supercharacter of the tensor product
$$S^1(V)^{\otimes m_1}\otimes S^2(V)^{\otimes m_2}\otimes \dots \otimes S^1(V^*)^{\otimes n_1}\otimes S^2(V^*)^{\otimes n_2}\dots ,$$ considered as a module over $A(n,n).$  When $i\ne 0$ the component $J(\mathfrak{g})_{i}$ is linearly generated by supercharacters of Kac modules.
The special case of $A(1,1)$ is considered separately in section 8.

In the $C(n)$ case due to Proposition \ref{cn} $J(\mathfrak
{g}) = \bigoplus_{a\in\mathbb C/\mathbb Z}J(\mathfrak{g})_{a},$ where again the components $J(\mathfrak{g})_{a}$ with $a \notin {\mathbb Z}$ are spanned over $\mathbb Z$ by the supercharacters of  Kac modules
 $K(\chi)$ with
$$
\chi= a\varepsilon+\sum_{j=1}^n\mu_{j}\delta_{j}, \mu_{1}\ge\mu_{2}\ge\dots\ge\mu_{n},\,\,\mu_{j}\in\Bbb Z_{\ge 0}.
$$
The zero component is the direct sum
$ J(\mathfrak{g})_{0}= J(\mathfrak{g})^{+}_{0}\oplus J(\mathfrak{g})^{-}_{0},$ where $J(\mathfrak{g})^{-}_{0}$ is spanned over $\mathbb Z$ by the supercharacters of  Kac modules $K(\chi)$ with
$$
\chi= \lambda\varepsilon+\sum_{j=1}^n\mu_{j}\delta_{j}, \mu_{1}\ge\mu_{2}\ge\dots\ge\mu_{n},\,\,\lambda\in\Bbb Z,\mu_{j}\in\Bbb Z_{\ge 0}
$$
and $J(\mathfrak{g})^{+}_{0}$ is generated over $\mathbb Z$ by $h_k,$ which are the supercharacters of  symmetric powers
of the standard representation. 

The proof of our main Theorem is now complete.

\section{Explicit description of the rings $J(\mathfrak{g})$}

In this section we describe explicitly the rings $J(\mathfrak{g})$ for all basic classical superalgebras except $A(1,1)$ case, which is to be considered separately in the next section. We start with the case of Lie superalgebra $\mathfrak{gl}(n,m),$ which will be used for the investigation of the $A(n,m)$ case.

\medskip

$\mathfrak{gl}(n,m)$

\medskip

%\noindent

In this case $\mathfrak{g}_{0}=\mathfrak{gl}(n)\oplus
\mathfrak{gl}(m)$ and $\mathfrak{g}_{1}=V_{1}\otimes V_{2}^*\oplus
V_{1}^*\otimes V_{2}$ where $V_{1}$ and $V_{2}$ are the identical
representations of $\mathfrak{gl}(n)$ and $\mathfrak{gl}(m))$
respectively. Let $\varepsilon_{1},\dots,\varepsilon_{n+m}$ be the
weights of the identical representation of  $\mathfrak{gl}(n,m)$.
Then the root system of $ \mathfrak{g}$ is expressed in terms of
linear functions $\varepsilon_{i},\ 1\le i\le n$ and
$\delta_{p}=\varepsilon_{p+n},\: 1\le p\le m$ as follows $$
R_{0}=\{\varepsilon_{i}-\varepsilon_{j}, \delta_{p}-\delta_{q} :
\: i\ne j\,:\: 1\le i,j\le n\ , p\ne q,\, 1\le p,q\le m\}, $$ $$
R_{1}=\{\pm(\varepsilon_{i}- \delta_{p}), \quad 1\le i\le n, \,
1\le p\le m\} = R_{iso}. $$ The invariant bilinear form is
determined by the relations $$
(\varepsilon_{i},\varepsilon_{i})=1,\:
(\varepsilon_{i},\varepsilon_{j})=0,\: i\ne j,\:
(\delta_{p},\delta_{q})=-1,\: (\delta_{p},\delta_{q})=0,\: p\ne
q,\: (\varepsilon_{i},\delta_{p})=0. $$ The Weyl group $W_0=
S_{n}\times S_{m}$ acts on the weights by separately
 permuting $\varepsilon_{i},\; i=1,\dots, n$ and
$\delta_{p},\; p=1,\dots, m$ .
Recall that the weight group of Lie algebra $\mathfrak{g}_{0}$ is defined as
 $$
 P_{0}=\{\lambda\in\mathfrak{h}^{*}\mid \frac{2(\lambda,\alpha)}{(\alpha,\alpha)}\in\mathbb Z
 \text{ for any } \alpha\in R_{0}\}.
  $$
In this case we have \begin{equation}
\label{pogl}
 P_{0}=\{\lambda\in\mathfrak{h}^{*}\mid
\lambda=\sum_{i=1}^m
\lambda_{i}\varepsilon_{i}+\sum_{p=1}^n\mu_{p}\delta_{p},\:
\lambda_{i}-\lambda_{j}\in \mathbb
Z\:\text{and}\:\mu_{p}-\mu_{q}\in \mathbb Z\} .
\end{equation}
 Choose the
following distinguished (in the sense of section 5) system of
simple roots $$
B=\{\varepsilon_{1}-\varepsilon_{2},\dots,\varepsilon_{n-1}-\varepsilon_{n},
\varepsilon_{n}- \delta_{1},\delta_{1}-
\delta_{2},\dots,\delta_{m-1}-\delta_{m}\}. $$ Note that the only
isotropic root is $\varepsilon_{n}- \delta_{1}.$ The weight
$\lambda$ is a highest weight for $\mathfrak{g}_{0}$ if
$\frac{2(\lambda,\alpha)}{(\alpha,\alpha)}\ge0$ for every non-isotropic
root $\alpha$ from $B.$

Let $x_i = e^{\varepsilon_i}, y_p = e^{\delta_p}$ be the elements
of the group ring of $\mathbb {Z}[P_0],$ which can be described as
the direct sum $\mathbb {Z}[P_0]= \bigoplus_{a,b\in\mathbb C/
\mathbb {Z}}\mathbb {Z}[P_0]_{a,b},$ where $$\mathbb
{Z}[P_0]_{a,b}=( x_{1}\dots x_{n})^{a}( y_{1}\dots y_{m})^{b}
\mathbb{Z}[x_{1}^{\pm1}\dots,x_{n}^{\pm1},y_{1}^{\pm1},\dots,y_{m}^{\pm
1} ]^{W_0}.$$ By definition the ring $J(\mathfrak{g})$ is the
subring $$J(\mathfrak{g})=\{f\in \mathbb {Z}[P_0]\mid
 y_{p}\frac{\partial f}{\partial y_{p}}+x_{i}\frac{\partial f}{\partial x_{i}} \in
 (y_{p}- x_{i} ), \quad p = 1, \dots, m, \quad i= 1, \dots, n\}.$$

Consider the rational function 
 $$
\chi(t)=\frac{\prod_{p=1}^m (1-y_{p}t)}{\prod_{i=1}^n(1-x_{i}t)}
$$ and expand it into Laurent series at zero and at infinity \footnote {The importance of considering the
Laurent series both at zero and infinity in this context was first
understood by Khudaverdian and Voronov \cite{KV}. They used this to write down some interesting relations in the Grothendieck ring of finite dimensional representations of $GL(m,n).$}  $$
\chi(t)=\sum_{k=0}^{\infty}h_{k}t^k=\sum_{k=n-m}^{\infty}h^{\infty}_{k}t^{-k}.
$$ 
Let us introduce
$$\Delta=\frac{y_{1}\dots y_{m}}{x_{1}\dots x_{n}},\quad
\Delta^*=\frac{x_{1}\dots x_{n}}{y_{1}\dots
y_{m}}=\Delta^{-1},\quad h_{k}^*=
h_{k}(x_{1}^{-1},\dots,x_{n}^{-1} ,y_{1}^{-1},\dots,y_{m}^{-1}).$$
It is easy to see that $h_{k}^{\infty}=\Delta h_{k+m-n}^*.$ 
We define also $h_k$ (and thus $h_k^{\infty}$) for all $k \in \mathbb Z$ by assuming that $h_k \equiv 0$ for negative $k.$

 \begin{proposition}\label{gl} The ring  $J(\mathfrak{g})$ for the Lie
 superalgebra $\mathfrak{gl}(n,m)$ is a direct sum
 $$
 J(\mathfrak{g})=\bigoplus_{a,b\in\mathbb C/\mathbb
 Z}J(\mathfrak{g})_{a,b},
 $$
where
 $$ J(\mathfrak{g})_{a,b}=( x_{1}\dots x_{n})^{a}( y_{1}\dots
y_{m})^{b}\prod_{i,p}(1-x_{i}/y_{p})\:\mathbb{Z}[x_{1}^{\pm1}\dots,x_{n}^{\pm1},y_{1}^{\pm1},\dots,y_{m}^{\pm
1} ]^{S_n\times S_m} $$ if $a+b\notin\mathbb Z;$

 $$ J(\mathfrak{g})_{a,b}= ( x_{1}\dots x_{n})^{a} (y_{1}\dots
y_{m})^{-a}J(\mathfrak{g})_{0,0} $$ if $a+b\in\mathbb Z, \quad
a\notin\mathbb Z$ and

$$
J(\mathfrak{g})_{0,0}=\{f\in\mathbb{Z}[x_{1}^{\pm1}\dots,x_{n}^{\pm1},y_{1}^{\pm1},\dots,y_{m}^{\pm
1} ]^{S_n\times S_m}\mid
 y_{p}\frac{\partial f}{\partial y_{p}}+x_{i}\frac{\partial f}{\partial x_{i}} \in
 (y_{p}- x_{i} ) \}.
$$
\end{proposition}

Proof easily follows from the definition of $J(\mathfrak{g}).$

\begin{proposition}\label{gl0}
The subring $J(\mathfrak{g})_{0,0}$ is generated over $\mathbb Z$
by $\Delta,\Delta^{*},\, h_{k},h^{*}_{k}, \,
k \in {\mathbb{N}}$ and can be interpreted as the Grothendieck ring of finite dimensional representations of algebraic supergroup $GL(n,m).$ 
\end{proposition}

\begin{proof}
We use the induction in $n+m$. When $n+m = 1$ it is obvious.
Assume that $n+m
> 1$. If $m=0$ or $n=0$ the statement follows from the
theory of symmetric functions \cite{Ma}. So we can assume that
$n>0$ and $m>0$. Consider a homomorphism $$ \tau :
J((\mathfrak{gl}(n,m))_{0,0} \longrightarrow
J(\mathfrak{gl}(n-1,m-1))_{0,0} $$ such that
$\tau(x_{n})=\tau(y_{m})= t$ and identical on others $x_i$ and
$y_p$. From the definition of $J((\mathfrak{gl}(n,m))$ it follows
that the image indeed belongs to $J(\mathfrak{gl}(n-1,m-1)).$ By
induction we may assume that $J(\mathfrak{gl}(n-1,m-1))_{0,0}$ is
generated by $\Delta,\:\Delta^{*}$ and $ h_{k},h_{k}^*$ for
$k=1,2,\dots$. We have $$
\tau(\Delta)(x_{1},\dots,x_{n-1}, t ,y_{1},\dots,y_{m-1}, t)=
\Delta(x_{1},\dots,x_{n-1},y_{1},\dots,y_{m-1}) $$ 
$$ \tau(
h_{k})(x_{1},\dots,x_{n-1}, t ,y_{1},\dots,y_{m-1}, t)=
 h_{k}(x_{1},\dots,x_{n-1},y_{1},\dots,y_{m-1})
$$ and the same for $\Delta^*$ and $h_{k}^{*},\: k=1,2,\dots$.
Therefore homomorphism $\tau$ is surjective. So now we need only
to prove that the kernel of $\tau$ is generated by
$\Delta,\Delta^{*}, h_{k},h_{k}^*$ for $k=1,2,\dots$.
 
Let
$a_0 = 1,\, a_{i}=(-1)^{i}\sigma_{i}(x),\, i = 1, \dots,n,$ where $\sigma_{i}$ are the elementary
symmetric polynomials in $x_{1},\dots, x_{n}$.  We have 
$$
\prod_{j=1}^m
(1-y_{j}t)=\chi(t)\sum_{i=0}^na_{i}t^i=\sum_{i=0}^na_{i}t^i\sum_{k\in\mathbb
Z}h_{k}t^k=\sum_{i=0}^na_{i}t^i\sum_{k\in\mathbb
Z}h^{\infty}_{-k}t^{k}.
$$
We see that
$$
\sum_{k\in\mathbb
Z}\left(\sum_{i=0}^nh_{k-i}a_{i}\right)t^k=\sum_{k\in\mathbb
Z}\left(\sum_{i=0}^nh^{\infty}_{-k+i}a_{i}\right)t^k,
$$
so we have the following infinite system of linear
equations  (see Khudaverdian and Voronov \cite{KV}):
$$
\sum_{i=0}^n(h_{k-i}-h^{\infty}_{i-k})a_{i}=0, \quad k\in\mathbb
{Z}.
$$

Introducing the elements $\tilde
h_{k}=h_{k}-h^{\infty}_{-k}$ we have
$$
\sum_{i=0}^n\tilde h_{k+n-i}a_{i}=0, \quad k=0,\:\pm1,\:\pm2,\:\dots$$
%Note that only for a
%finite number of $k$ both $h_{k}$ and $h^{\infty}_{-k}$ are
%non-zero.
Considering this as a linear system for the unknown $a_1, \dots, a_n$ with given $a_0=1$   we have by  Cramer's rule for any pairwise  different $k_{1},\dots,k_{n}$
$$
\left|\begin{array}{cccc}
 \tilde h_{k_{1}}& \tilde h_{k_{1}+1}& \ldots &\tilde h_{k_{1}+n-1}\\
\tilde h_{k_{2}}& \tilde h_{k_{2}+1}& \ldots &\tilde h_{k_{2}+n-1}\\
\vdots&\vdots&\ddots&\vdots\\
\tilde h_{k_{n}}& \tilde h_{k_{n}+1}& \ldots &\tilde h_{k_{n}+n-1}\\
 \end{array}\right|a_{n}=(-1)^n
\left|\begin{array}{cccc}
 \tilde h_{k_{1}+1}& \tilde h_{k_{1}+2}& \ldots &\tilde h_{k_{1}+n}\\
\tilde h_{k_{2}+1}& \tilde h_{k_{2}+2}& \ldots &\tilde h_{k_{2}+n}\\
\vdots&\vdots&\ddots&\vdots\\
\tilde h_{k_{n}+1}& \tilde h_{k_{n}+2}& \ldots &\tilde h_{k_{n}+n}\\
 \end{array}\right|
$$
and more  generally  for any integer $l$
\begin{equation}
\label{det}
\left|\begin{array}{cccc}
 \tilde h_{k_{1}}& \tilde h_{k_{1}+1}& \ldots &\tilde h_{k_{1}+n-1}\\
\tilde h_{k_{2}}& \tilde h_{k_{2}+1}& \ldots &\tilde h_{k_{2}+n-1}\\
\vdots&\vdots&\ddots&\vdots\\
\tilde h_{k_{n}}& \tilde h_{k_{n}+1}& \ldots &\tilde h_{k_{n}+n-1}\\
 \end{array}\right|a^{l}_{n}=(-1)^{nl}
\left|\begin{array}{cccc}
 \tilde h_{k_{1}+l}& \tilde h_{k_{1}+l+1}& \ldots &\tilde h_{k_{1}+n+l-1}\\
\tilde h_{k_{2}+l}& \tilde h_{k_{2}+l+1}& \ldots &\tilde h_{k_{2}+n+l-1}\\
\vdots&\vdots&\ddots&\vdots\\
\tilde h_{k_{n}+l}& \tilde h_{k_{n}+l+1}& \ldots &\tilde h_{k_{n}+n+l-1}\\
 \end{array}\right|
\end{equation}
Any element from kernel of $\tau$ has a form
$$f = R(x,y) g(x,y),\quad
g\in\mathbb{Z}[x_{1}^{\pm1}\dots,x_{n}^{\pm1},y_{1}^{\pm1},\dots,y_{m}^{\pm
1} ]^{S_n\times S_m},$$ where $$R(x,y)= \prod_{i=1}^{n}\prod_{p=1}^m \left(1-\frac{y_p}{x_i}\right).$$ 
Let $s_{\lambda}(x),s_{\mu}(y)$ be the Schur functions corresponding to the sequences of non-increasing  integers $\lambda = ( \lambda_{1}\ge\dots\ge \lambda_{n})$, $\mu= ( \mu_{1}\ge\dots\ge \mu_{m})$ (see \cite{Ma}). It is easy to see that the products $s_{\lambda}(x) s_{\mu}(y)$
give a basis in $\mathbb{Z}[x_{1}^{\pm1}\dots,x_{n}^{\pm1},y_{1}^{\pm1},\dots,y_{m}^{\pm
1} ]^{S_n\times S_m}.$

Thus we need to show that $f_{\lambda,\mu}=
s_{\lambda}(x)s_{\mu}(y)R(x,y)
$ can be expressed in terms of $ h_{k},\: h_{k}^*, \:\Delta,\:\Delta^*.$  Multiplying $f_{\lambda,\mu}$ by an appropriate power of $\Delta$ we can
assume that
$
f_{\lambda,\mu}=a_{n}^{l}s_{\lambda}(x)s_{\mu}(y)R(x,y),
$
where $l$ is an integer and $\lambda,\mu$
are partitions (i.e. $\lambda_n$ and $\mu_m$ are non-negative) such that $\lambda_{n}\ge m$. But in this case we can use the well-known formula (see e.g. \cite{Ma}, I.3, Example 23)
$$
s_{\lambda}(x)s_{\mu}(y)R(x,y)=
\left|\begin{array}{cccc}
  h_{\lambda_{1}}&  h_{\lambda_{1}+1}& \ldots & h_{\lambda_{1}+p+n-1}\\
 h_{\lambda_{2}-1}&  h_{\lambda_{2}}& \ldots & h_{\lambda_{2}+p+n-2}\\
\vdots&\vdots&\ddots&\vdots\\
 h_{\lambda_{n}-n+1}&  h_{\lambda_{n}-n+2}& \ldots & h_{\lambda_{n}+p}\\
  h_{\mu_{1}^{\prime}-n}&  h_{\mu_{1}^{\prime}-n+1}& \ldots & h_{\mu_{1}^{\prime}+p-1}\\
 \vdots&\vdots&\ddots&\vdots\\
 h_{\mu^{\prime}_{p}-p-n+1}& h_{\mu^{\prime}_{p}-p-n+2}& \ldots &h_{\mu^{\prime}_{p}}\\
 \end{array}\right|
 $$
where $\mu^{\prime}_{1},\dots,\mu^{\prime}_{p}$ be the partition conjugated to $\mu_{1},\dots,\mu_{m}$.
Since $\lambda_{n}\ge m$ for any $h_{k}$ from the first $n$ rows we have $h_{k}=\tilde h_{k}$.
Let us multiply this equality by $a^{l}_{n}$ and then expand the determinant
with respect to the first $n$ rows  by Laplace's rule. Using (\ref{det}) we get
$$
f_{\lambda,\mu}=\left|\begin{array}{cccc}
  \tilde h_{\lambda_{1}+l}&  \tilde h_{\lambda_{1}+l+1}& \ldots &\tilde h_{\lambda_{1}+l+p+n-1}\\
 \tilde h_{\lambda_{2}+l-1}&  \tilde h_{\lambda_{2}+l}& \ldots & \tilde h_{\lambda_{2}+l+p+n-2}\\
\vdots&\vdots&\ddots&\vdots\\
  h_{\mu^{\prime}_{p}-p-n+1}&  h_{\mu^{\prime}_{p}-p-n+2}& \ldots &
  h_{\mu^{\prime}_{p}}\\
 \end{array}\right|.
 $$ 
 Thus we have shown that $J(\mathfrak{g})_{0,0}$ is generated 
by $\Delta,\Delta^{*},\, h_{k},h^{*}_{k}.$ Since all these elements are the supercharacters of some representations of the algebraic supergroup $GL(n,m)$ (see e.g. \cite{BK}) we see that $J(\mathfrak{g})_{0,0}$ is a subring of the 
Grothendieck ring of this supergroup. Other elements of  $J(\mathfrak{g})$ can not be extended already to the algebraic subgroup $GL(n)\times GL(m)$, so $J(\mathfrak{g})_{0,0}$ coincides with the
Grothendieck ring of $GL(n,m).$
\end{proof}

Now we are going through the list of basic classical Lie superalgebras.

\medskip

$A(n-1,m-1)$

\medskip

%\noindent

%We have the following
\begin{proposition}\label{amnneq}
The ring  $J(\mathfrak{g})$ for the Lie
 superalgebra $\mathfrak{sl}(n,m)$ with $(n,m)\ne (2,2)$
is a direct sum
 $$
 J(\mathfrak{g})=\bigoplus_{a\in\mathbb C/\mathbb Z}J(\mathfrak{g})_{a},
 $$
\begin{equation}
\label{pr1}
J(\mathfrak{g})_{a}=\{ f\in( x_{1}\dots x_{n})^{a}\prod_{i,p}(1-x_{i}/y_{p})\mathbb{Z}[x^{\pm1}, y^{\pm1}]^{S_n\times S_m}_0
\end{equation}
if $a\notin\mathbb Z$ and
\begin{equation}
\label{jo}
J(\mathfrak{g})_{0}=\{f\in\mathbb{Z}[x^{\pm1},y^{\pm 1}]^{S_n\times S_m}_0\mid
 y_{j}\frac{\partial f}{\partial y_{j}}+x_{i}\frac{\partial f}{\partial x_{i}} \in
 (y_{j}- x_{i} ) \},
\end{equation}
where $
\mathbb{Z}[x^{\pm1},y^{\pm 1}]^{S_n\times S_m}_0$ is the quotient of the ring 
 $\mathbb{Z}[x_{1}^{\pm1}\dots,x_{n}^{\pm1},y_{1}^{\pm1},\dots,y_{m}^{\pm 1}
]^{S_n\times S_m}$ by the ideal generated by $x_1\dots x_n - y_1\dots y_m.$

The subring $J(\mathfrak{g})_{0}$ is generated over $\mathbb Z$ by $h_{k},h^{*}_{k},\, k\in \mathbb{N}$
and  can be interpreted as the Grothendieck ring of finite dimensional representations of algebraic supergroup $SL(n,m).$ 
\end{proposition}

The first part easily follows from Proposition \ref{gl}, the description of $J(\mathfrak{g})_{0}$ is based on Proposition \ref{gl0}.
The case $m=n$ is special.

\medskip

$A(n-1,n-1)=\mathfrak{psl}(n,n),\, n >2.$

\medskip

%\noindent

The root system of $A(n-1,n-1)$ is
$$
R_{0}=\{\tilde\varepsilon_{i}-\tilde\varepsilon_{j}, \tilde\delta_{p}-\tilde\delta_{q} :
\: i\ne j\,\: 1\le i,j\le n\ , p\ne q,\, 1\le p,q\le n\} $$ $$
R_{1}=\{\pm(\tilde\varepsilon_{i}- \tilde\delta_{p}), \quad 1\le i\le n, \,
1\le p\le n\}  = R_{iso}$$
where
$$
\tilde\varepsilon_{1}+\dots+\tilde\varepsilon_{n}=0,\:\tilde\delta_{1}+\dots+\tilde\delta_{n}=0.
$$
These weights are related to the weights of $\frak{sl}(n,n)$ by the formulas
$$\tilde\varepsilon_{i}= \varepsilon_{i}-\frac{1}{n} \sum_{j=1}^n \varepsilon_{j},\quad \tilde\delta_{i}= \delta_{i}-\frac{1}{n} \sum_{j=1}^n \delta_{j}, \, i=1,\dots, n.$$
The bilinear form is defined by the relations
$$
(\tilde\varepsilon_{i},\tilde\varepsilon_{i})=1-1/n,\:
(\tilde\varepsilon_{i},\tilde\varepsilon_{j})=-1/n,\: i\ne j,$$
$$
(\tilde\delta_{p},\tilde\delta_{p})=-1+1/n,\: (\tilde\delta_{p},\tilde\delta_{q})=1/n,\:
p\ne q,\: (\tilde\varepsilon_{i},\tilde\delta_{p})=0.
$$
The Weyl group $W_0= S_{n}\times S_{n}$ acts on the weights
by  permuting separately
$\tilde\varepsilon_{i},\; i=1,\dots,n$ and
$\tilde\delta_{p},\; p=1,\dots,n$ . 
A distinguished system of simple roots can be chosen as
$$
B=\{\tilde\varepsilon_{1}-\tilde\varepsilon_{2},\dots,\tilde\varepsilon_{n-1}-\tilde\varepsilon_{n},
\tilde\varepsilon_{n}- \tilde\delta_{1},\tilde\delta_{1}-
\tilde\delta_{2},\dots,\tilde\delta_{n-1}-\tilde\delta_{n}\}.
$$
The weight  lattice of the Lie
algebra $\mathfrak{g}_{0}$ is
$$
P_{0}=\{\sum_{i=1}^{n-1}
\lambda_{i}\tilde\varepsilon_{i}+\sum_{p=1}^{n-1}\mu_{p}\tilde\delta_{p} \mid
\lambda_{i}, \mu_p \in \mathbb
Z\}.
$$
\begin{proposition} \label{ann}
The ring $J(\frak g)$ for Lie superalgebra $\frak g=  \frak{psl}(n,n)$ with $n > 2$ is a direct sum
$$
J(\frak g)=\bigoplus_{i=0}^{n-1} J(\frak g)_{i}
$$
where for $i\ne0$
$$
J(\frak g)_{i}=\{ f=( x_{1}\dots x_{n})^{\frac{i}{n}}\prod^n_{j,p}(1-x_{j}/y_{p})g,\, g \in \mathbb{Z}[x^{\pm1},y^{\pm 1}]^{S_n\times S_n}_0, \, \deg g= -i\},
$$
and $J(\frak g)_{0}$ is the subring of (\ref{jo}) with $m=n$, consisting of elements of degree $0.$
 
 The  ring  $J(\frak g)_{0}$ is  linearly generated by the products
$$
h_{1}^{m_{1}}h_{2}^{m_{2}}\dots( h^*_{1})^{n_{1}}( h^*_{2})^{n_{2}}\dots
$$
such that $m_{1}+2m_{2}+\dots= n_{1}+2n_{2}+\dots$ and can be interpreted as  the Grothendieck ring of finite dimensional representations of the algebraic supergroup $PSL(n,n)$.
\end{proposition} 

\begin{proof}
From the definition of the ring $J(A(n-1,n-1))$ it follows that this ring can be identified with the subring in $J(\frak{sl}(n,n))$ consisting of the linear combinations of 
$$e^{\lambda_{1}\varepsilon_{1}+\dots+\lambda_{n}\varepsilon_{n}+\mu_{1}\delta_{1}+\dots+\mu_{n}\delta_{n}}
$$
such that $\lambda_{1}+\dots+\lambda_{n}+\mu_{1}+\dots+\mu_{n}=0$. This subring can be also characterised as the ring of invariants with respect to the automorphism 
$$
\theta_{t}(x_{i})=tx_{i},\,\,\theta_{t}(y_{i})=ty_{i}
$$
Now the proposition easy follows from these formulas and proposition \ref{amnneq}.
\end{proof}

\medskip
$C(n)=\mathfrak{osp}(2,2n)$

\medskip

%\noindent.

In this case $\mathfrak{g}_{0}=\mathfrak{so}(2)\oplus sp(2n)$ and $\mathfrak{g}_{1}=V_{1}\otimes
V_{2},$ where $V_{1}$ and $V_{2}$ are  the identical representations of $so(2)$ and $sp(2n)$ respectively. 

Let $\varepsilon_{1},\dots,\varepsilon_{n+1}$ be the weights of the identical representation of  $C(n)$ and define
$\varepsilon=\varepsilon_{1},\,
\delta_{j}=\varepsilon_{j+1},\: 1\le j\le n.$
The root system is
 $$
R_{0}=\{\pm\delta_{i}\pm\delta_{j},\:\pm2\delta_{i}, \: i\ne j,\: 1\le i,j\le n\}
$$
$$
R_{1}=\{\pm\varepsilon\pm\delta_{j},\:\pm\delta_{j}\},\quad
R_{iso}=\{\pm\varepsilon\pm\delta_{j},\}
$$
with the bilinear form
$$
(\varepsilon,\varepsilon)=1,\:
(\delta_{i},\delta_{i})=-1,\: (\delta_{i},\delta_{j})=0,\: i\ne j,\: (\varepsilon,\delta_{k})=0
$$
The Weyl group $W_0$ is the semi-direct product of $S_{n}$ and ${Z}_{2}^n$. It acts on the weights by  permuting and changing the signs of $ \delta_{j},j=1,\dots,n$.
As a distinguished system of simple roots we select
$$
B=\{\varepsilon-\delta_{1},\delta_{1}-\delta_{2},\dots,\delta_{n-1}-\delta_{n},2\delta_{n}\}.
$$

The weight group has the form
$$
P_{0}=\{ \nu=\lambda\varepsilon+\sum_{j=1}^n\mu_{j}\delta_{j},\:
\lambda \in \mathbb C,\,\mu_{j}\in \mathbb Z\}.
$$
Let $e^{\varepsilon}=x,\: e^{\delta_{j}}=y_{j},\: u=x+x^{-1},\: v_{j}=y_{j}+y_{j}^{-1},\: j=1,\dots, n$.
 Consider the Taylor expansion at zero
of the following rational function 
$$
\chi(t)=\frac{\prod_{j=1}^m (1-y_{j}t) (1-y_{j}^{-1}t)}{(1-xt) (1-x^{-1}t)}=\sum_{k=0}^{\infty}h_{k}t^k.$$

\begin{proposition} \label{cn}
The ring  $J(\mathfrak{g})$ for the Lie
 superalgebra $C(n)$
is a direct sum
$$
 J(\mathfrak{g})=\bigoplus_{a\in\mathbb C/\mathbb Z}J(\mathfrak{g})_{a},
$$
where
$$
J(\mathfrak{g})_{a}= x^{a}\prod_{j=1}^n(1-x/y_{j})(1-xy_{j})\:\mathbb{Z}[x^{\pm1},y_{1}^{\pm1},\dots,y_{n}^{\pm 1}
]^{W_0}
$$
if $a\notin\mathbb Z$ and
 $$
 J(\mathfrak{g})_{0}=\{f\in\mathbb{Z}[x^{\pm1}, y_{1}^{\pm1},\dots,y_{n}^{\pm 1}
]^{W_0}\mid
 y_{j}\frac{\partial f}{\partial y_{j}}+x\frac{\partial f}{\partial x} \in
 (y_{j}- x),\, j=1,\dots, n \}.
 $$
More explicitly,
$
 J(\mathfrak{g})_{0}= J(\mathfrak{g})^{+}_{0}\oplus J(\mathfrak{g})^{-}_{0},
$
where
$$
 J(\mathfrak{g})^{-}_{0}=\{ f=x\prod_{j=1}^n(u-v_{j})g\mid\: g\in\mathbb{Z}[u,v_{1},\dots,v_{n}]^{S_{n}}\},
$$
$$
 J(\mathfrak{g})^{+}_{0}=\{f\in\mathbb{Z}[u,v_{1},\dots,v_{n}]^{ S_{n}}\mid
 u\frac{\partial f}{\partial u}+ v_{j}\frac{\partial f}{\partial v_{j}}\in (u-v_{j}),\, j=1,\dots n\}.
$$

The subring $J(\mathfrak{g})^{+}_{0} $ is generated over $\mathbb Z$ by $ h_{k}, \, k \in \mathbb N$ and can be interpreted as the Grothendieck ring of finite dimensional representations of the algebraic supergroup $OSP(2,2n).$ 
\end{proposition}

\begin{proof}
The first claim is obvious.
To prove the second one note that $x^2-xu+1=0$. Therefore any element $f$ from
$ J(\mathfrak{g})_{0}$ can be uniquely written in the form $f_{0}+xf_{1}$, where $f_0,\,f_1\in \mathbb{Z}[u,v_{1},\dots,v_{n}]^{S_{n}}$. Condition $ y_{j}\frac{\partial f}{\partial y_{j}}+x\frac{\partial f}{\partial x} \in (y_j-x)$ means that after substitution $y_{j}=x$ the polynomial $f=f_{0}+xf_{1}$ does not depend on $x$. Because of the symmetry $y_j \rightarrow y_j^{-1}$ the same must be true for $f_{0}+x^{-1}f_{1}.$ This means that $f_{1}$ is zero after substitution $y_{j}=x$, which implies the claim.

The fact that $J(\mathfrak{g})^{+}_{0} $ is generated by $ h_{k}$ follows from the theory of supersymmetric functions \cite{Ma}. 

Since $h_k$ are the supercharacters of the symmetric powers of the standard representation all elements of $J(\mathfrak{g})^{+}_{0}$ give rise to representations of the supergroup $OSP(2,2n).$ The elements of $J(\mathfrak{g})^{-}_{0}$ can not be extended already to the subgroup $O(2).$
\end{proof}

$B(m,n)=\mathfrak{osp}(2m+1,2n)$

\medskip

%\noindent
Here $\mathfrak{g}_{0}=so(2m+1)\oplus sp(2n)$ and
$\mathfrak{g}_{1}=V_{1}\otimes V_{2}$ where $V_{1}$ and $V_{2}$ are  the
identical representations of $so(2m+1)$ and $sp(2n)$ respectively. Let
$\pm\varepsilon_{1},\dots,\pm\varepsilon_{m},\, \pm\delta_1, \dots, \pm\delta_n$ be the non-zero weights of the
identical representation of  $B(m,n).$
Then the root system of $B(m,n)$ is 
$$R_{0}=\{\pm\varepsilon_{i}\pm\varepsilon_{j},\:\pm\varepsilon_{i}, \, \pm\delta_{p}\pm\delta_{q},\:\pm2\delta_{p}, \,
\: i\ne j,\,\: 1\le i,j\le m\ , p\ne q,\, 1\le p,q\le n\} 
$$
$$
R_{1}=\{\pm\varepsilon_{i}\pm\delta_{p},\:\pm\delta_{p}\},\quad
R_{iso}=\{\pm\varepsilon_{i}\pm\delta_{p}\}.
$$
The invariant bilinear form is
$$
(\varepsilon_{i},\varepsilon_{i})=1,\: (\varepsilon_{i},\varepsilon_{j})=0,\: i\ne j,\:
(\delta_{p},\delta_{p})=-1,\: (\delta_{p},\delta_{q})=0,\: p\ne q,\: (\varepsilon_{i},\delta_{p})=0.
$$
The Weyl group $W_0=\left(S_{n}\ltimes\mathbb{Z}_{2}^n\right)\times\left( S_{m}\ltimes\mathbb{Z}_{2}^m\right)$ acts on the weights by separately  permuting $\varepsilon_{i},\; j=1,\dots,m$ and
$\delta_{p},\; p=1,\dots,n$ and changing their signs.
The weight  lattice of the Lie algebra $\mathfrak{g}_{0}$ is
$$
P_{0}=\{ \nu =\sum_{i=1}^m\lambda_{i}\varepsilon_{i}+\sum_{p=1}^n\mu_{p}\delta_{p},\:\lambda_{i}\in \mathbb Z\:\text{or}\:\lambda_{i}\in \mathbb Z+\frac{1}{2}\, \text{for all $i$},\;\:
\mu_{p}\in \mathbb Z\}.
$$
A distinguished system of simple roots can be chosen as
$$
B=\{\delta_{1}-\delta_{2},\dots,\delta_{n-1}-\delta_{n},\delta_{n}-\varepsilon_{1},\varepsilon_{1}-\varepsilon_{2},\dots,\varepsilon_{m-1}-\varepsilon_{m}, \varepsilon_{m}\}.
$$
The weight $\lambda$ is a highest weight of $\mathfrak{g}_{0}$ if $\frac{2(\lambda,\alpha)}{(\alpha,\alpha)}\ge0$ for any simple root of $\mathfrak{g}_{0}$
 $$
 \alpha\in\{\delta_{1}-\delta_{2},\dots,\delta_{n-1}-\delta_{n}, 2\delta_{n},\varepsilon_{1}-\varepsilon_{2},\dots,\varepsilon_{m-1}-\varepsilon_{m}, \varepsilon_{m}\}
 $$
 Introduce the variables $x_{i}=e^{\varepsilon_{i}},\:x^{1/2}_{i}=e^{\varepsilon_{i}/2},\, u_{i}=x_{i}+x_{i}^{-1},\, i=1,\dots,m$ and $ y_{p}=e^{\delta_{p}},\:  v_{p}=y_{p}+y_{p}^{-1},\: p=1,\dots, n$.
 Consider the Taylor series at zero of the following function
$$
\chi(t)=\frac{\prod_{p=1}^n (1-y_{p}t) (1-y_{p}^{-1}t)}{(1-t)\prod_{i=1}^m(1-x_{i}t) (1-x_{i}^{-1}t)}=\sum_{k=0}^{\infty}h_{k}(x,y)t^k
$$

\begin{proposition} \label{bmn}
The ring  $J(\mathfrak{g})$ of Lie
 superalgebra of type $B(m,n)$
is a direct sum
$$
 J(\mathfrak{g})=J(\mathfrak{g})_{0} \oplus J(\mathfrak{g})_{1/2},
$$
where
$$
J(\mathfrak{g})_{1/2}= \prod_{i=1}^m(x^{1/2}_{i}+x^{-1/2}_{i})
\prod_{i,p}(u_{i}-v_{p})g\mid g\in\mathbb{Z}[u_{1},\dots, u_{m},v_{1},\dots,v_{n}]^{S_{m}\times S_{n}}
$$
and
$$
J(\mathfrak{g})_{0}= \{f\in\mathbb{Z}[u_{1},\dots, u_{m},v_{1},\dots,v_{n}]^{S_{m}\times S_{n}}\mid
 u_{i}\frac{\partial f}{\partial u_{i}}+ v_{p}\frac{\partial f}{\partial v_{p}}\in (u_{i}-v_{p})\}.
$$
The subring  $J(\mathfrak{g})_{0} $ is generated over $\mathbb Z$ by $ h_{k}(x,y), \, k\in \mathbb Z$
and can be interpreted as the Grothendieck ring of finite dimensional representations of the algebraic supergroup $OSP(2m+1,2n).$
\end{proposition}

\begin{proof}
The decomposition  $J(\mathfrak{g})=J(\mathfrak{g})_{0} \oplus J(\mathfrak{g})_{1/2} $ reflects the fact that all $\lambda_i$ in the weight lattice $P_0$ are either all integer or half-integers. Consider $f\in J(\mathfrak{g})$ and suppose first that all the corresponding $\lambda_i$ are half integer. Write $f$ as a Laurent polynomial with respect to $x_{1},y_{1}$
 $$
 f=\sum c_{i,j}x_{1}^iy_{1}^j,
 $$
 where the coefficients $c_{i,j}$ depend on the remaining variables. The condition $ x_{1}\frac{\partial f}{\partial x_{1}}+ y_{1}\frac{\partial f}{\partial y_{1}}\in (x_{1}-y_{1})$
 means  that $\sum(i+j) c_{i,j}=0$. Since $i$ is not an integer but $j$ does we conclude that $\sum c_{i,j}=0$. This means that $f$ is divisible by $(x_{1}-y_{1})$ and hence by the symmetry by $\prod_{i,p}(u_{i}-v_{p})$. The factor $\prod_{i=1}^m(x^{1/2}_{i}+x^{-1/2}_{i})$ is due to the Weyl group symmetry of $B(m).$
The last part is similar to the previous case..
\end{proof}

$D(m,n)=\mathfrak{osp}(2m,2n),\, m>1$

\medskip

%\noindent
In this case $\mathfrak{g}_{0}=so(2m)\oplus sp(2n)$ and $\mathfrak{g}_{1}=V_{1}\otimes
V_{2},$ where $V_{1}$ and $V_{2}$ are the identical representations of $so(2m)$ and $sp(2n)$ respectively. Let $\pm\varepsilon_{1},\dots,\pm\varepsilon_{m},\, \pm\delta_1, \dots, \pm\delta_n$ be the weights of the identical representation of  $D(m,n)$.
The root system is 
$$R_{0}=\{\pm\varepsilon_{i}\pm\varepsilon_{j}, \, \pm\delta_{p}\pm\delta_{q},\:\pm2\delta_{p}, \,
\: i\ne j,\,\: 1\le i,j\le m\ , p\ne q,\, 1\le p,q\le n\} 
$$
$$
R_{1}=\{\pm\varepsilon_{i}\pm\delta_{p}\} =
R_{iso}.
$$
The bilinear form is defined by the relations
$$
(\varepsilon_{i},\varepsilon_{i})=1,\: (\varepsilon_{i},\varepsilon_{j})=0,\: i\ne j,\:
(\delta_{p},\delta_{p})=-1,\: (\delta_{p},\delta_{q})=0,\: p\ne q,\: (\varepsilon_{i},\delta_{p})=0.
$$
The Weyl group $W=\left(S_{m}\ltimes\mathbb{Z}_{2}^{m-1}\right)\times
\left( S_{n}\ltimes\mathbb{Z}_{2}^{n}\right)$ acts on the
weights by separately  permuting 
$\varepsilon_{i},\; i=1,\dots,m$ and $\delta_{p},\; p=1,\dots,n$ and changing their signs 
such that the total change of signs of $\varepsilon_i$ is even.

The weight  lattice of Lie algebra $\mathfrak{g}_{0}$ is the same as in the previous case:
$$
P_{0}=\{ \nu =\sum_{i=1}^m\lambda_{i}\varepsilon_{i}+\sum_{p=1}^n\mu_{p}\delta_{p},\:\lambda_{i}\in \mathbb Z\:\text{or}\:\lambda_{i}\in \mathbb Z+\frac{1}{2}\, \text{for all $i$},\;\:
\mu_{p}\in \mathbb Z\}.
$$
A distinguished system of simple roots is
$$
B=\{\delta_{1}-\delta_{2},\dots,\delta_{n-1}-\delta_{n},\delta_{n}-\varepsilon_{1},\varepsilon_{1}-\varepsilon_{2},\dots,\varepsilon_{m-1}-\varepsilon_{m}, \varepsilon_{m-1}+\varepsilon_{m}\}
$$
The weight $\nu$ is a highest weight for $\mathfrak{g}_{0}$ if $\frac{2(\nu,\alpha)}{(\alpha,\alpha)}\ge0$ for all simple roots of $\mathfrak{g}_{0}$
 $$
 \alpha\in\{\delta_{1}-\delta_{2},\dots,\delta_{n-1}-\delta_{n}, 2\delta_{n},\varepsilon_{1}-\varepsilon_{2},\dots,\varepsilon_{m-1}-\varepsilon_{m}, \varepsilon_{m-1}+\varepsilon_{m}\},
 $$
which is equivalent to
 $$
 \mu_{1}\ge\dots\ge\mu_{n}\ge0,\:\lambda_{1}\ge\dots\ge\lambda_{m-1}\ge|\lambda_{m}|.
 $$
 Introduce the variables $x_{i}=e^{\varepsilon_{i}},\:x^{1/2}_{i}=e^{\varepsilon_{i}/2},\, u_{i}=x_{i}+x_{i}^{-1},\, i=1,\dots,m$ and $ y_{p}=e^{\delta_{p}},\:  v_{p}=y_{p}+y_{p}^{-1},\: p=1,\dots, n$ and consider the following Taylor series
$$
\chi(t)=\frac{\prod_{p=1}^n (1-y_{p}t) (1-y_{j}^{-1}t)}{\prod_{i=1}^m(1-x_{i}t) (1-x_{i}^{-1}t)}=\sum_{k=0}^{\infty}h_{k}(x,y)t^k.
$$
We will need also the following invariant of the Weyl group $D(m)$
$$
\omega=\sum x_{1}^{\pm1}\dots x_{m}^{\pm1},
$$
where the sum is over all possible combinations of $\pm 1$ with even sum.

\begin{proposition} \label{dmn}
The ring  $J(\mathfrak{g})$ of Lie
 superalgebra of type $D(m,n)$
is a direct sum
$$
 J(\mathfrak{g})=J(\mathfrak{g})_{0} \oplus J(\mathfrak{g})_{1/2},
$$
where
$$
J(\mathfrak{g})_{1/2}=\{
\prod_{i,p}(u_{i}-v_{p})\left((x_{1}\dots x_{m})^{1/2}\mathbb{Z}[x_{1}^{\pm1},\dots, x_{m}^{\pm1},y_{1}^{\pm1},\dots, y_{n}^{\pm1}]\right)^{W_0}\}
$$
and
$$
  J(\mathfrak{g})_{0}=\{f\in\mathbb{Z}[x_{1}^{\pm1},\dots,x_{m}^{\pm1} ,y_{1}^{\pm1},\dots,y_{n}^{\pm 1}
]^{W_0}\mid
 y_{p}\frac{\partial f}{\partial y_{p}}+x_{i}\frac{\partial f}{\partial x_{i}} \in
 (y_{p}- x_{i} )\}
 $$
More explicitly,
$$
 J(\mathfrak{g})_{0}= J(\mathfrak{g})^{+}_{0}\oplus J(\mathfrak{g})^{-}_{0},
$$
where
$$
 J(\mathfrak{g})^{-}_{0}=\{\omega\prod_{i,p}(u_{i}-v_{p})\mathbb{Z}[u_{1}\dots,u_{m},v_{1},\dots,v_{n}]^{S_{m}\times S_{n}}\},
$$
$$
J(\mathfrak{g})^{+}_{0}=\{f\in\mathbb{Z}[u_{1},\dots,u_{m},v_{1},\dots v_{n}]^{ S_{m}\times S_{n}}\mid
 u_i\frac{\partial f}{\partial u_i}+ v_{p}\frac{\partial f}{\partial v_{p}}\in (u_i-v_{p})\}
$$
The subring  $J(\mathfrak{g})^{+}_{0} $ is generated over $\mathbb Z$ by $ h_{k}(x,y),\, k\in\mathbb Z$ and can be interpreted as the Grothendieck ring of finite dimensional representations of algebraic supergroup $OSP(2m,2n)$.
\end{proposition}

\begin{proof}
The proof of the first claim is similar to the previous case. Let us explain the decomposition of $J(\mathfrak{g})_{0}.$
 It is well known (see e.g. \cite {FH}) that any element from $\mathbb{Z}[x_{1}^{\pm1},\dots, x_{m}^{\pm1},y_{1}^{\pm1},\dots, y_{n}^{\pm1}]^{W_0}$ can be written uniquely in the form  $f_{0}+\omega f_{1}$,
where $f_{0}, f_{1}\in \mathbb{Z}[u_{1},\dots,u_{m},v_{1},\dots v_{n}]^{ S_{m}\times S_{n}}.$ The condition $ y_{1}\frac{\partial f}{\partial y_{1}}+x_{1}\frac{\partial f}{\partial x_{1}}\in(x_{1}-y_{1})$ means that after the substitution $y_{1}=x_1=t$ the polynomial $f$ does not depend on $t$. Because of the symmetry $y_1 \rightarrow y_1^{-1}$ the same must be true for $f_{0}+\tau(\omega)f_{1},$ where the transformation $\tau$ maps $x_{1}$ to $x_{1}^{-1}$ and leaves the remaining variables invariant.  
Therefore $(\omega-\tau(\omega)) f_{1}$ does not depend on $t$ after the substitution $y_{1}=x_1=t$. This implies that $f_{1}$ is zero after this substitution, which explains the form of  $J(\mathfrak{g})^{-}_{0}$.
The last part is standard by now. 
\end{proof}

$G(3)$

\medskip

%\noindent

In this case $\mathfrak{g}_{0}=G(2)\oplus\mathfrak{sl}(2)$ and
$\mathfrak{g}_{1}=U \otimes V$, where $U$
is the first fundamental representation of $G(2)$ (see \cite{OV} or \cite{Bou1})
and $V$ is the identity representation of  $\mathfrak{sl}(2)$. Let
 $\pm\varepsilon_{i}, i=1,2,3,\,\varepsilon_{1}+\varepsilon_{2}+\varepsilon_{3}=0$ be the non-zero weights of $U$
and $\pm\delta$ be the weights of identity representation of $\mathfrak{sl}(2)$.
Then the root system of $G(3)$ is
$$
R_{0}=\{\varepsilon_{i}-\varepsilon_{j},\:\pm\varepsilon_{i},\:\pm2\delta\},\quad
R_{1}=\{\pm\varepsilon_{i}\pm\delta,\:\pm\delta\},\quad
R_{iso}=\{\pm\varepsilon_{i}\pm\delta\}
$$
with the bilinear form defined by
$$
(\varepsilon_{i},\varepsilon_{i})=2\quad(\varepsilon_{i},\varepsilon_{j})=-1,\: i\ne j,
\quad (\delta,\delta)=-2.
$$
The Weyl group $W_0=
D_6 \times\mathbb{Z}_{2}$, where
$D_6$ is the dihedral group of order 12 acting on $\varepsilon_{i}$ by
permutations and simultaneously changing their signs, while $\mathbb{Z}_{2}$ is acting by changing the sign of $\delta$.
 The weight lattice of the Lie algebra $\mathfrak{g}_{0}$ can be written as
$$
P_{0}=\{\nu=\lambda_{1}\varepsilon_{1}+\lambda_{2}\varepsilon_{2}+\mu\delta
,\:\lambda_{1},\lambda_{2},\mu\in\mathbb Z\}.
$$
A distinguished system of simple roots is
 $$
 B=\{\varepsilon_{3}+\delta,\: \varepsilon_{1},\:\varepsilon_{2}-\varepsilon_{1}\}.
 $$
 The weight $\lambda$ is a highest weight for $\mathfrak{g}_{0}$ if $(\lambda,\alpha)\ge0$ for $\alpha\in\{
 \varepsilon_{1},\:\varepsilon_{2}-\varepsilon_{1},\delta\}$, which is equivalent to the following conditions
 $\lambda_{1} \ge \lambda_{2}-\lambda_{1}\ge0$. Let $x_{1}=e^{\varepsilon_{1}},\: x_{2}=e^{\varepsilon_{2}}, \: y=e^{\delta}.$  By definition the ring  $J(\mathfrak{g})$ is
$$
 J(\mathfrak{g})=
\{f\in\mathbb{Z}[x_{1}^{\pm1},x_{2}^{\pm1}, y^{\pm1}]^{W_0}\mid
 x_{1}\frac{\partial f}{\partial x_1} +x_{2}\frac{\partial f}{\partial x_{2}}+2y\frac{\partial f}{\partial y}\in
 (y- x_{1} x_{2}) \},$$
where the action of $W_0$ is generated by the permutation of $x_1$ and $x_2$ and
by the transformations $x_1\rightarrow (x_1x_2)^{-1}, \, x_2\rightarrow x_2$ and $x_1\rightarrow x_1^{-1},\, x_2\rightarrow x_2^{-1}.$
Let
$
 u_{1}=x_{1}+x_{1}^{-1},\,  u_{2}=x_{2}+x_{2}^{-1},\,  u_{3}=x_{1}x_2+x_{1}^{-1}x_2^{-1},\, v=y+y^{-1}$ and introduce
 $$
 w=v^{2}-v(u_{1}+u_{2}+u_{3}+1)+u_{1}u_{2}+u_{1}u_{3}+u_{2}u_{3},
 $$
which (up to additional constant 1) is the supercharacter of the adjoint representation, and hence belongs to $J(\mathfrak{g}).$

 \begin{proposition} \label{g12}
The ring  $J(\mathfrak{g})$ of Lie
 superalgebra of type $G(3)$ can be described as
$$
J(\mathfrak{g})=\{ f=g(w)+(v-u_{1})(v-u_{2})(v-u_{3})h\mid h\in\mathbb{Z}[u_{1},u_{2},u_{3},v]^{S_{3}},\, g\in \mathbb{Z}[w]\}.
 $$
 \end{proposition}

\begin{proof}
 It is not difficult to verify that
$$
\mathbb{Z}[x_{1}^{\pm1},x_{2}^{\pm1},y^{\pm1}]^{\mathbb{Z}_{2}\times\mathbb{Z}_{2}}=
\mathbb{Z}[u_{1},u_{2},u_{3},v],
$$
where the generators of $\mathbb{Z}_{2}\times\mathbb{Z}_{2}$ are acting by changing $y\rightarrow y^{-1}$ and $x_1 \rightarrow x_1^{-1},\,x_2 \rightarrow x_2^{-1}.$ 
For any $ f\in  J(\mathfrak{g})$ consider $q=f(x_{1}, x_{2}, x_{1}x_{2}),$ then we have
$$
x_{1}\frac{\partial f}{\partial x_1} +x_{2}\frac{\partial f}{\partial x_{2}}+2y\frac{\partial f}{\partial y}=
x_{1}\frac{\partial q}{\partial x_1} +x_{2}\frac{\partial q}{\partial x_{2}}=0$$
when $ y=x_{1}x_{2}.$ This means that $q$ has degree $0.$ Since $q$ is also invariant under the transformation $x_{1}\rightarrow x_{1}^{-1},\: x_{2}\rightarrow x_{2}^{-1}$ there exists a polynomial of one variable  $g$ such that $q = g(\frac{x_{1}}{x_{2}}+\frac{x_{2}}{x_{1}})$. But it is easy to check that when $y=x_{2}x_{3}$ then $w=\frac{x_{1}}{x_{2}}+\frac{x_{2}}{x_{1}}.$ Therefore  the  difference
$f-g(w)$ is divisible by  $(y- x_{1} x_{2})$ and by the symmetry it is also divisible by
$$
(y- x_{1} x_{2})(y- x_{1}^{-1} x_{2}^{-1})(y- x_{1} )(y- x_{1}^{-1})(y- x_{2})(y-x_{2}^{-1})=
y^3(v-u_{1})(v-u_{2})(v-u_{3}).
$$
A simple check shows that any polynomial of the form $(v-u_{1})(v-u_{2})(v-u_{3})h,\,h\in\mathbb{Z}[u_{1},u_{2},u_{3},v]^{S_{3}}$ belongs to $ J(\mathfrak{g}).$ 
\end{proof}

$F(4)$

\medskip

%\noindent
In this case $\mathfrak{g}_{0}=B_{3}\oplus\mathfrak{sl}(2)$
and $\mathfrak{g}_{1}=U\otimes V$, where $U$
is the spin representation of $B_{3}$ (see \cite{OV} or \cite{Bou1})
and $V$ is the identity representation of  $\mathfrak{sl}(2)$. Let $\pm\varepsilon_{1},\pm\varepsilon_{2}, \pm\varepsilon_{3}$ are the non-zero weights of the identity
representation of $B_{3}$ and $\pm\frac{1}{2}\delta$ be the weights of identity representation of $\mathfrak{sl}(2).$ The root
system of $\mathfrak{g}$ is
$$
R_{0}=\{\pm\varepsilon_{i}\pm\varepsilon_{j},\:\pm\varepsilon_{i}\:\pm\delta\},\quad
R_{1}=\{\frac{1}{2}(\pm\varepsilon_{1}\pm\varepsilon_{2}\pm\varepsilon_{3}\pm\delta)\} =
R_{iso}
$$
with the bilinear form defined by
$$
(\varepsilon_{i},\varepsilon_{i})=1\quad(\varepsilon_{i},\varepsilon_{j})=0,\: i\ne j,
\quad (\delta,\delta)=-3.
$$
The Weyl group $W_0$ is $ (S_{3}\ltimes\mathbb{Z}^3_{2})\times\mathbb{Z}_{2}$, where
$S_{3}\ltimes\mathbb{Z}^3_{2}$ acts on $\varepsilon_{i}$'s by
permutations and changing their signs while the second factor $\mathbb{Z}_{2}$
changes the sign of $\delta$.
 The weight lattice  of the Lie algebra $\mathfrak{g}_{0}$ is
$$
P_{0}=\{ \nu=\lambda_{1}\varepsilon_{1}+\lambda_{2}\varepsilon_{2}+\lambda_{3}\varepsilon_{3}+\mu\delta,\:
\lambda_{i}\in\mathbb Z \:\text{or}\:\lambda_{i}\in\mathbb Z+1/2,\:2\mu\in\mathbb Z\}
$$
As a distinguished system of simple root we choose
 $$
 B=\{\frac{1}{2}(\delta-\varepsilon_{1}-\varepsilon_{2}-\varepsilon_{3}), \varepsilon_{1}-\varepsilon_{2},\varepsilon_{2}-\varepsilon_{3},\varepsilon_{3}\}.
 $$
 The weight $\lambda$ is a highest weight for $\mathfrak{g}_{0}$ if $\frac{(\nu,\alpha)}{(\alpha,\alpha)}\ge0$ for $\alpha\in\{
 \varepsilon_{1}-\varepsilon_{2},\:\varepsilon_{2}-\varepsilon_{3},\varepsilon_{3},\;\delta\}$, which is equivalent to the following conditions:
 $\lambda_{1}\ge \lambda_{2}\ge\lambda_{3}\ge0,\: \mu\ge0$.  Let $x_{1}=e^{\frac{1}{2}\varepsilon_{1}},\: x_{2}=e^{\frac{1}{2}\varepsilon_{2}},\: x_{3}=e^{\frac{1}{2}\varepsilon_{3}},\: y=e^{\frac{1}{2}\delta}, \,
 u_{i}=x_{i}+x_{i}^{-1} (i=1,2,3),\: v=y+y^{-1}.
 $
By definition the ring  $J(\mathfrak{g})$ consists of polynomials
 $f\in\mathbb{Z}[x_{1}^{\pm2},x_{2}^{\pm2},x_{3}^{\pm2},
( x_{1}x_{2}x_{3})^{\pm 1},y^{\pm1}]^{W_0}$ such that 
$$
 3y\frac{\partial f}{\partial y}+x_{1}\frac{\partial f}{\partial x_{1}} +x_{2}\frac{\partial f}{\partial x_{2}}+x_{3}\frac{\partial f}{\partial x_{3}}\in
 (y- x_{1}x_{2}x_{3}).
 $$
 Introduce 
 $$
  Q=(v-x_{1}x_{2}x_{3}-x^{-1}_{1}x_{2}^{-1}x_{3}^{-1})\prod_{i=1}^3
 \left(v-\frac{x_{1}x_{2}x_{3}}{x_{i}^2}-\frac{x_{i}^2}{x_{1}x_{2}x_{3}}\right)
 $$
and
 $$
 w_{k}= \sum_{i\ne j} \frac{x_{i}^{2k}}{x_{j}^{2k}}+\sum_{i=1}^3(x_{i}^{2k}+x_{i}^{-2k})+y^{2k}+y^{-2k}-
 (y^{k}+y^{-k})\prod_{i=1}^3(x_{i}^k+x_{i}^{-k}), \,\,k=1,2.$$
 It is easy to check that $Q h$ belongs to the ring $J(\mathfrak{g})$ for any polynomial 
$h$ from $\mathbb{Z}[x_{1}^{\pm2},x_{2}^{\pm2},x_{3}^{\pm2},( x_{1}x_{2}x_{3})^{\pm 1},y^{\pm1}]^{W_0}.$ The element $w_1$ up to a constant is the supercharacter of the adjoint representation, $w_2$ can be expressed as a linear combination of the supercharacters
of the tensor square of the adjoint representation and its second symmetric power, so both of them also belong to the ring.

\begin{proposition}
The ring  $J(\mathfrak{g})$ of Lie
 superalgebra of type $F(4)$ can be described as
$$
 J(\mathfrak{g}) = \{f=g(w_{1},w_{2})+Qh \mid
 h\in\mathbb{Z}[x_{1}^{\pm2},x_{2}^{\pm2},x_{3}^{\pm2},( x_{1}x_{2}x_{3})^{\pm 1},y^{\pm1}]^{W_0},\, g\in\mathbb{Z}[w_1, w_2]\}.$$
\end{proposition}

\begin{proof}

Let $ f\in J(\mathfrak{g})$  and consider $q=f(x_{1},x_{2},x_{3}, x_{1}x_{2} x_{3}).$ We have
$$
3y\frac{\partial f}{\partial y}+x_{1}\frac{\partial f}{\partial x_{1}}+x_{2}\frac{\partial f}{\partial x_{2}} +x_{3}\frac{\partial f}{\partial x_{3}}=x_{1}\frac{\partial g}{\partial x_{1}}+
x_{2}\frac{\partial g}{\partial x_{2}} +x_{3}\frac{\partial g}{\partial x_{3}}=0$$ when $ y=x_{1}x_{2}x_{3}.
$
This means as before that $q$ has degree $0$ and therefore it is a Laurent polynomial in $x_{1}^{2},x_{2}^{2},x_{3}^{2}.$ Since $q$ is invariant under the transformations $x_{i}\rightarrow x_{i}^{-1}$ and the permutation group $S_{3}$ there exists a polynomial $g$ of two variables such that $q = g(u_1, u_2)$, where
$$
u_1 = \sum_{i\ne j}\frac{x_{i}^2}{x_{j}^2},\:
u_2 = \sum_{i\ne j}\frac{x_{i}^4}{x_{j}^4}.
$$
But it is easy to check that when $y=x_{1}x_{2}x_{3}$ we have 
$
w_{1}=u_1,\, w_{2}= u_2.
$
Therefore  the difference
$f-g(w_{1}, w_{2})$ is divisible by  $(y- x_{1}x_{2} x_{3})$ and by symmetry is divisible by $Q$.
\end{proof}

$D(2,1,\alpha)$

\medskip

In this case $\mathfrak{g}_{0}=\mathfrak{sl}(2)\oplus\mathfrak{sl}(2)
\oplus\mathfrak{sl}(2),\,\mathfrak{g}_{1}=V_{1}\otimes V_{2}\otimes V_{3}$,
where $V_{i}$ are the identity representations of the corresponding $\mathfrak{sl}(2)$.
Let $\pm\varepsilon_{1}, \pm\varepsilon_{2}, \pm\varepsilon_{3}$ be their weights. The root system of $\mathfrak{g}$ is 
$$
R_{0}=\{\pm2\varepsilon_{1},\pm2\varepsilon_{2},\pm2\varepsilon_{3}\}\quad
R_{1}=\{\pm\varepsilon_{1}\pm\varepsilon_{2}\pm\varepsilon_{3}\},\quad
R_{iso}=R_{1}.
$$
The bilinear form depends on the parameter $\alpha$:
$$
(\varepsilon_{1},\varepsilon_{1})=-1-\alpha,\:(\varepsilon_{2},\varepsilon_{2})=1,\:
(\varepsilon_{3},\varepsilon_{3})=\alpha
$$
The Weyl group $W_0= \mathbb{Z}^{3}_{2}$ acts on the $\varepsilon_{i}$'s by
 changing their signs.
 The weight lattice  of the Lie algebra $\mathfrak{g}_{0}$ is
 $$
P_{0}=\{\lambda=\lambda_{1}\varepsilon_{1}+\lambda_{2}
\varepsilon_{2}+\lambda_{3}\varepsilon_{3},\:\lambda_{1},\lambda_{2},\lambda_{3}
\in\mathbb Z\}.
$$
 Choose the following distinguished system of simple roots
 $$
 B=\{\varepsilon_{1}+\varepsilon_{2}+\varepsilon_{3},\: -2\varepsilon_{2},\:-2\varepsilon_{3}
 \},
 $$
then the highest weights $\lambda$ satisfy following conditions:
$\lambda_{1}\ge0,\:\lambda_{2}\le0,\:\lambda_{3}\le0$.

Let  $x_{i}=e^{\varepsilon_{i}},\: u_{i}=x_{i}+x_{i}^{-1},\: i=1,2,3$. 
By definition  we have $$ J(\mathfrak{g})=\{f\in\mathbb{Z}[x_{1}^{\pm1},x_{2}^{\pm1},x_{3}^{\pm1}]^{\mathbb{Z}^{3}_{2}}\mid
 (1+\alpha)x_{1}\frac{\partial f}{\partial x_{1}}+x_{2}\frac{\partial f}{\partial x_{2}} +
 \alpha x_{3}\frac{\partial f}{\partial x_{3}}
 \in
 (x_{1}-x_{2}x_{3}) \}.
 $$
Introduce
 $$
 Q=(x_{1}-x_{2}x_{3})(x_{2}-x_{1}x_{3})(x_{3}-x_{1}x_{2})(1-x_{1}x_{2}x_{3})x_{1}^{-2}x^{-2}_{2}
 x^{-2}_{3}=
u_{1}^2+u_{2}^2+u_{3}^2-u_{1} u_{2} u_{3}-4,
 $$
which up to a constant is the supercharacter of the adjoint representation.
For the rational non-zero values of the parameter $\alpha= p/q, \, p \in\mathbb{Z}, \,q\in\mathbb{N}$ we will need the additional element
   $$
 w_{\alpha}=(x_{1}+x_{1}^{-1}-
 x_{2}x_{3}-x_{2}^{-1}x_{3}^{-1})
 \frac{(x_{2}^p-x_{2}^{-p})(x_{3}^q-x_{3}^{-q})}{(x_{2}-x_{2}^{-1})(x_{3}-x_{3}^{-1})}+x_{2}^{p}x_{3}^{-q}+ x_{2}^{-p}x_{3}^{q},
 $$
which also belongs to $J(\mathfrak{g})$ as one can check directly.

\begin{proposition} \label{d21alpha}
 If $\alpha$ is not rational then the ring $J(\mathfrak{g})$ of the Lie superalgebra $D(2,1,\alpha)$ can be described as follows
$$
J(\mathfrak{g}) = \{f=c+Q h \mid c \in \mathbb{Z},\, h\in\mathbb{Z}[u_{1},u_{2},u_{3}]\}.
 $$
 If $\alpha =p/q$ is rational  then 
 $$
J(\mathfrak{g}) =\{ f=g(w_{\alpha})+Q h\mid h\in\mathbb{Z}[u_{1},u_{2},u_{3}], \, g \in 
\mathbb{Z}[w]\}. $$
\end{proposition}

\begin{proof}
First note that
$
\mathbb{Z}[x_{1}^{\pm1},x_{2}^{\pm1},x_{3}^{\pm1}]^{\mathbb{Z}^{3}_{2}}=
\mathbb{Z}[u_{1},u_{2},u_{3}].
$
Take $ f\in J(\mathfrak{g})$  and consider the function $\phi(x_2, x_3)=f(x_{2}x_{3}, x_{2}, x_{3}),$ then
$$
(1+\alpha)x_{1}\frac{\partial f}{\partial x_{1}}+x_{2}\frac{\partial f}{\partial x_{2}} +
 \alpha x_{3}\frac{\partial f}{\partial x_{3}}=
x_{2}\frac{\partial \phi}{\partial x_{2}} +\alpha x_{3}\frac{\partial \phi}{\partial x_{3}}=0$$
when $x_{1}=x_{2}x_{3}.$
If $\alpha$ is irrational then $\phi$ must be a constant. If $\alpha= p/q$ is rational
 then $\phi = g(x_{2}^px_{3}^{-
q}+x_{2}^{-p}x_{3}^{ q})$  for some polynomial $g \in \mathbb{Z}[w]$ since it is invariant under the transformation $x_{2}\rightarrow
x_{2}^{-1},\: x_{3}\rightarrow x_{3}^{-1}$.
But when $x_{1}=x_{2}x_{3}$ the element
$w_{\alpha}=x_{2}^px_{3}^{- q}+x_{2}^{-p}x_{3}^{ q}$.
Therefore  the difference $f - g(w_{\alpha})$ is divisible by $(x_{1}-
x_{2} x_{3})$ and by symmetry by $Q$.
\end{proof}

\section{Special case $A(1,1)$}

This case is special because the isotropic roots have multiplicity 2.
The definition of the ring $J(\mathfrak{g})$  should be modified in this case as follows:
\begin{equation}
\label{Jdef} J(\mathfrak{g})=\{ f\in\mathbb{Z}[ P]^{W_0}:  \,
D_{\alpha}f\in ((e^{\alpha}-1)^2) \quad \text{for any isotropic
root } \alpha \}
\end{equation}
where $((e^{\alpha}-1)^2)$ denotes the principal ideal in $
\mathbb{Z}[ P]$ generated by $(e^{\alpha}-1)^2.$
We would like to note that the property (\ref{Jdef}) can be rewritten as
$$D_{\alpha} \frac{1}{(e^{\alpha}-1)} D_{\alpha}f \in (e^{\alpha}-1),$$
which is a natural analogue of the condition proposed for the quantum Calogero-Moser
systems by Chalykh and one of the authors in \cite{CV1990}.

\begin{thm} \label{a11}
The Grothendieck ring $K(\mathfrak{g})$ of
finite dimensional representations of  Lie superalgebra
$\mathfrak{g}=A(1,1)=\mathfrak{psl}(2,2)$ is isomorphic to the ring $J(\mathfrak
{g}).$ The isomorphism is given by the supercharacter map
$Sch: K(\mathfrak{g}) \rightarrow J(\mathfrak {g}).$
\end{thm}

Now we are going to prove this result.
We have in this case $\mathfrak{g}_{0}=\mathfrak{sl}(2)\oplus\mathfrak{sl}(2),\,\mathfrak{g}_{1}=V_{1}\otimes V_{2}\oplus V_{1}\otimes V_{2},$
where $V_{1}, V_{2}$ are the identity representations of the corresponding $\mathfrak{sl}(2)$.
Let $\{\varepsilon,-\varepsilon,\delta,-\delta\}$
be the corresponding weights. The roots of $A(1,1)$
are
$$
R_{0}=\{2\varepsilon,\:-2\varepsilon, \: 2\delta,\:-2\delta\}
$$
$$
R_{iso}=R_{1}=\{\varepsilon+\delta,\:
\varepsilon-\delta,\:-\varepsilon+\delta,\: -\varepsilon-\delta\}
$$
and the invariant bilinear form is
$$
(\varepsilon,\varepsilon)=1, \: (\delta,\delta)=-1,\:
(\varepsilon,\delta)=0.
$$
The important fact is that the multiplicity of any isotropic root
equals to two. Note that in this case $R$ is not a generalized root system in the sense
of the definition given in section 3: one can check that the third property is not satisfied.
However it is the case if we use a more general definition proposed by Serganova \cite{Serga} and used in our previous work \cite{SV}.

The weight lattice of
the Lie algebra $\mathfrak{g}_{0}$ is
$$
P_{0}=\{ \nu=
\lambda\varepsilon+\mu\delta,\: \lambda,\mu\in \mathbb Z\}
$$
The Weyl group $W_0= \mathbb{Z}_{2}\times \mathbb{Z}_{2}$ is acting
on the weights by changing the signs of $\varepsilon$ and $\delta$.  
A distinguished system of simple roots is
$$
B=\{\varepsilon-\delta,2\delta\}.
$$
The weight $\nu=\lambda\varepsilon+\mu\delta$ is a highest weight for
$\mathfrak{g}_{0}$ if $\lambda,\mu\ge 0$. 

The following result generalizes the proposition \ref{J} to the case when multiplicities of the isotropic roots are equal to 2.
 
\begin{proposition}\label{J2}
 Let  $\mathfrak{g}$
be the solvable Lie superalgebra such that
$\mathfrak{g}_{0}=\mathfrak{h}$ is a commutative finite
dimensional Lie algebra, $\mathfrak{g}_{1}=Span(X_{1}, X_{2},
Y_{1}, Y_{2})$ and the following relations hold
$$ [h,X_{i}]=\alpha(h)X_{i},\:
[h,Y_{i}]=-\alpha(h)Y_{i},\: [ Y_{i},Y_{j}]=[ X_{i},X_{j}]=0,\: [
X_{i},Y_{j}]=\delta_{i,j}H,\: i,j=1,2
$$
where $H \in \mathfrak{h}$ and $\alpha\ne0$ is a linear form on
$\mathfrak{h}$ such that $ \alpha(H)=0.$ Then the Grothendieck
ring of $\mathfrak{g}$ is isomorphic to
\begin{equation}
\label{Jprop} J(\mathfrak{g})=\{ f=\sum c_{\lambda}
e^{\lambda}\mid \lambda\in\mathfrak{h}^*,\quad D_{H}f\in
((e^{\alpha}-1)^2))\}.
\end{equation}
The
isomorphism is given by the supercharacter map
$
Sch: [V] \longrightarrow sch \,V.
$
\end {proposition}

\begin{proof}
Every irreducible finite-dimensional $\mathfrak{g}$-module $V$  has unique (up to a
multiple) vector $v$ such that $X_{1} v=X_{2} v = 0, \, h v =
\lambda(h)v$ for some linear form $\lambda$ on $\mathfrak{h}.$
This establishes a bijection between the  irreducible
$\mathfrak{g}$-modules and the elements of $\mathfrak{h}^*.$

There are two types of such modules, depending on whether
$\lambda(H)=0$ or not.  In the first case the module
$V=V({\lambda})$ is one-dimensional and its supercharacter is
$e^{\lambda}$. If $\lambda(H)\ne0$ then the corresponding module
$V({\lambda})$ is four-dimensional with the  supercharacter
$sch(V) = e^{\lambda}(1-e^{-\alpha})^2.$ In both cases the supercharacters belong to the ring $J(\mathfrak{g})$. Thus we have proved that
the image of $Sch\, (K(\mathfrak{g}))$ is
contained in $J(\mathfrak{g})$.

Conversely, let $f=\sum c_{\lambda}e^{\lambda}$ belong to
$J(\mathfrak{g})$. By subtracting a suitable linear combination of
supercharacters of the one-dimensional modules $V(\lambda)$ we can
assume that $\lambda(H)\ne0$ for all $\lambda$ from $f$.
The
condition $D_{H}f \in
((e^{\alpha}-1)^2))$ implies that $D_{H}f \in (e^{\alpha}-1).$ Using the same arguments
as in the proof of Proposition \ref{J}  we
deduce that $f$ itself belongs to the ideal generated by $(e^{\alpha}-1)$.
This means that $f =
(e^{\alpha}-1)h$ for some $h \in \mathbb{Z}[\mathfrak{h}^*]$.
It is easy to see that
$D_{H}f\in ((e^{\alpha}-1)^2)$ is equivalent to the condition
$D_{H}h\in (e^{\alpha}-1)$. Therefore as a before $h \in (e^{\alpha}-1),$ so
$f \in ((e^{\alpha}-1)^2).$  From the proof of the first part we conclude that  $f$ is a linear combination of the
supercharacters of the irreducible $\mathfrak{g}$-modules.
\end{proof}

Let
$x=e^{\varepsilon},\:y=e^{\delta},\: u=x+x^{-1},\:v=y+y^{-1}$.
By definition we have
$$
J(\mathfrak g)=\{f\in \mathbb Z[x^{\pm1},y^{\pm1}]^{\mathbb{Z}_{2}\times \mathbb{Z}_{2}}\mid
x\frac{\partial f}{\partial x}+y\frac{\partial f}{\partial y}\in
((x-y)^2)\}.
$$
\begin{proposition}
The ring  $J(\mathfrak{g})$ of Lie
 superalgebra of type $A(1,1)$ can be described as
 $$
J(\mathfrak g)=\{f=c+(u-v)^2g(u,v) \mid c \in \mathbb Z,\, g\in \mathbb Z[u,v]\}.
$$
The subring  $J(\mathfrak{g})^+$ of polynomials of even degree in $J(\mathfrak g)$ can be interpreted as the Grothendieck ring of finite dimensional representations of algebraic supergroup $PSL(2,2)$.
\end{proposition}

\begin{proof} The isomorphism
$$
\mathbb Z[x^{\pm1},y^{\pm1}]^{\mathbb{Z}_{2}\times \mathbb{Z}_{2}}=\mathbb Z[u,v]
$$
is standard. Take any $f\in J(\mathfrak g)$. We can write $f$ in the form
$
f=c+(u-v)q(v)+(u-v)^2g(u,v)
$
for some $c \in \mathbb Z,\, q \in \mathbb Z[v],\, g\in \mathbb Z[u,v].$ From the identity
$
u-v=(x-y)(1-1/xy)
$
it follows that $(u-v)^2g(u,v)\in J(\mathfrak g)$. Therefore $(u-v)q(v)\in
J(\mathfrak g)$. But it is easy to verify that in this case $q$ must be zero.

Let us prove the statement about $PSL(2,2)$. From the isomorphism $A(1,1)=\frak{psl}(2,2)$  we have the natural imbedding 
$$
J(A(1,1))=K(A(1,1))\longrightarrow K(\frak{sl}(2,2)) = J(\frak{gl}(2,2))/I,
$$
such that
$$
u\rightarrow \left(\frac{x_1}{x_2}\right)^{\frac12}+\left(\frac{x_2}{x_1}\right)^{\frac12},\,\,\,
v\rightarrow \left(\frac{y_1}{y_2}\right)^{\frac12}+\left(\frac{y_2}{y_1}\right)^{\frac12}
$$
and $I$ is the ideal generated by $1-e^{a(\varepsilon_{1}+\varepsilon_{2}-\delta_{1}-\delta_{2})},\,a\in\Bbb C$.
In the same way   as in the proposition \ref{ann} one can prove that the ring  $K(PSL(2,2))$  can be identified with the subring in $J(\frak{sl}(2,2))$ linearly generated by $h_{i_1}\dots h_{i_s}h^*_{j_1}\dots h^*_{j_r}$ such that $i_{1}+\dots+i_{s}=j_{1}+\dots+j_{r}$. But it is not difficult to verify that this subring coincides with the image of $J^+(\mathfrak g)$. Proposition is proved.
\end{proof}

Now the Theorem \ref{a11} follows from the fact that any polynomial of the form
$(u-v)^2\chi_k(u)\chi_l(v)$, where $\chi_k(u),\,\chi_l(v)$ are the characters of
the irreducible $A(1)$ modules with the highest weights $k$ and $l$, is the supercharacter of  a Kac module over $A(1,1)$.

\section{Super Weyl groupoid}

In this section we associate to any generalized root system (in Serganova's sense)
$R \subset V$ a certain groupoid $\mathfrak{W}= \mathfrak{W}(R),$ which we will call super Weyl groupoid.\footnote {We should note that the possibility of a groupoid version of the Weyl group for Lie superalgebras was contemplated by Serganova \cite{Serga3}, but she had a different picture in mind (see \cite{Serga4}).
Recently Heckenberger and Yamane \cite{Hecken2} introduced a groupoid related to basic classical Lie superalgebras motivated by Serganova's work and notion of the Weyl groupoid for Nichols algebras \cite{Hecken1}. Our super Weyl groupoid has no direct relations with this.}\,The corresponding Grothendieck ring can be interpreted as the invariant ring of a natural action of this groupoid. 

For a nice introduction to the theory of groupoids, including some history,  we refer to  the surveys by Brown  \cite{B} and Weinstein  \cite{W}. Recall that a {\it groupoid} can be defined as a small category with all morphisms being invertible. The set of objects is denoted as $\mathfrak B$ and called the {\it base} while the set of morphisms is denoted as $\mathfrak G.$  We will follow the common tradition to use the same notation $\mathfrak G$ for the groupoid itself.

If the base $\mathfrak B$ consists of one element $\mathfrak G$ has a group structure. More generally, for any $x \in \mathfrak B$ one can associate an {\it isotropy group} $\mathfrak G_x$ consisting of all morphisms $g\in \mathfrak G$ from $x$ into itself. For any groupoid we have a natural equivalence relation on the base $\mathfrak B$, when $x \sim y$ if there exists a morphism $g\in \mathfrak G$ from $x$ to $y.$ One can think therefore of groupoids as generalisations of both groups and the equivalence relations. 
In fact, any finite groupoid is a disjoint union of its subgroupoids called {\it components}, corresponding to the equivalence classes called {\it orbits}. Each such component up to an isomorphism is uniquely determined by the orbit and its isotropy group (see \cite{B}).

A standard example of groupoid comes from the action of a group $\Gamma$ on a set $X:$ the base $\mathfrak B=X$ and set $\mathfrak G$ of morphisms from $x$ to $y$ consists of the elements $\gamma \in \Gamma$ such that $\gamma(x) = y.$ 
%We will denote the corresponding groupoid as $(X, \Gamma).$ 

One can generalize this example in the following way.  Let $\mathfrak G$ be a groupoid and the group $\Gamma$ is acting on it by the automorphisms of the corresponding category. In particular, $\Gamma$ acts on the base $\mathfrak B$ of $\mathfrak G$
(for convenience, on the right).
Then one can define a {\it semi-direct product groupoid}
$\Gamma \ltimes \mathfrak G$ with the same base $\mathfrak B$ and the morphisms from $x$ to $y$ being pairs
$(\gamma,\,f),\, \gamma \in \Gamma, f \in \mathfrak G$ such that $f : \gamma(x) \rightarrow y.$ The composition is defined in a natural way: $(\gamma_1,\,f_1) \circ (\gamma_2,\,f_2) = (\gamma_1 \gamma_2,\,  \gamma_2(f_1) \circ f_2).$

Now we are ready to define the super Weyl groupoid $\mathfrak{W}(R)$ corresponding to generalized root system $R$.  Recall that the reflections with respect to the non-isotropic roots generate a finite group denoted $W_0.$ 

Consider first the following groupoid $\mathfrak T_{iso}$ with the base $R_{iso},$ which is the set of all the isotropic roots in $R.$ The set of morphisms from $\alpha \rightarrow \beta$ is non-empty if and only if $\beta = \pm \alpha$ in which case it consists of just one element. We will denote the corresponding morphism $\alpha \rightarrow -\alpha$ as $\tau_{\alpha}, \alpha \in R_{iso}.$
The group $W_0$ is acting on $\mathfrak T_{iso}$ in a natural way: $\alpha \rightarrow w(\alpha),\,
\tau_{\alpha} \rightarrow \tau_{w(\alpha)}.$ 
We define now the {\it super Weyl groupoid} $$\mathfrak{W}(R) = W_0 \coprod  W_0 \ltimes \mathfrak T_{iso}$$ as a disjoint union of the group $W_0$ considered as a groupoid with a single point base $[W_0]$ and the semi-direct product groupoid $W_0 \ltimes \mathfrak T_{iso}$ with the base $R_{iso}.$ Note that the disjoint union is a well defined operation on the groupoids.

There is a natural action of the groupoid $\mathfrak{W}(R)$ on the ambient space $V$ of generalized root system $R$ in the following sense.

For any set $X$ one can define the following groupoid $\mathfrak S(X)$, whose base consists of all possible subsets $Y \subset X$ and the morphisms are all possible bijections between them. By the {\it action of a groupoid $\mathfrak G$ on a set} $X$ we will mean the homomorphism  of $\mathfrak G$ into $\mathfrak S(X)$ (which is a functor between the corresponding categories). In case if $X=V$ is a vector space and $Y \subset X$ are the affine subspaces with morphisms being affine bijections, we will talk about {\it affine action}.

Let $X= V$ and define the following affine action $\pi$ of the super Weyl groupoid
$\mathfrak{W}(R)$ on it. 
The base point $[W_0]$ maps to the whole space $V$, while the base element corresponding to an isotropic root $\alpha$ maps to the hyperplane $\Pi_{\alpha}$ defined by the equation $(\alpha,x)=0.$
The elements of the group $W_0$ are acting in a natural way and the element $\tau_{\alpha}$ acts as a shift 
$$\tau_{\alpha}(x) = x + \alpha, \, x \in \Pi_{\alpha}.$$
Note that since $\alpha$ is isotropic $x + \alpha$ also belongs to $\Pi_{\alpha}.$ One can easily check that this indeed defines an affine action
of $\mathfrak{W}(R)$ on $V.$ 

A version of this action can be seen in the definition of the algebra $\Lambda_{R,B}$ of quantum integrals of the deformed Calogero-Moser systems introduced in our paper \cite{SV}: in that case the element $\tau_{\alpha}$ acts as a shift between two different affine hyperplanes (see formula (7) in \cite{SV}). The above defintion of the super Weyl groupoid was mainly motivated by this action.

The following reformulation of our main theorem shows that the super Weyl groupoid may be considered as a substitute of the Weyl group in the theory of Lie superalgebras. 

Let $V = \mathfrak h^*$ be the dual space to a Cartan subalgebra $\mathfrak h$
of a basic classical Lie superalgebra $\mathfrak g$ with generalized root system $R.$
Using the invariant bilinear form we can identify $V$ and $V^* = \mathfrak h$
and consider the elements of the group ring $\mathbb{Z}[\mathfrak {h}^*]$ as functions on $V.$
A function $f$ on $V$ is {\it invariant under the action of groupoid} $\mathfrak{W}$
if for any $g \in \mathfrak{W}$ we have $f(g(x)) = f(x)$ for all $x$ from the definition domain of the action map of $g.$

Let $P_0 \subset \mathfrak{h}^*$ be the abelian group of weights of $\mathfrak{g}_{0}$  and ${\mathbb{Z}}[P_0]$ be the corresponding integral group ring. It is easy to see 
that the ring $J(\mathfrak{g})$ is nothing else but the invariant elements of ${\mathbb{Z}}[P_0]$ 
invariant under the action $\pi$ of the super Weyl groupoid described above. Thus our main result can be reformulated as follows.

\begin{thm} 
The Grothendieck ring $K(\mathfrak{g})$ of the finite dimensional representations of a basic classical Lie superalgebra $\mathfrak{g}$ except $A(1,1)$ is isomorphic to the ring  ${\mathbb{Z}}[P_0]^{\mathfrak{W}}$ of invariants of the super Weyl groupoid $\mathfrak{W}$ under the action defined above. 
\end{thm}

\section{Concluding remarks}

Thus we have now a description of the Grothendieck rings of finite-dimensional representations for all basic classical Lie superalgebras. The fact that the corresponding rings can be described by simple algebraic conditions seems to be remarkable.
We believe that these rings as well as the corresponding super Weyl groupoids will play an important role in the representation theory. 

An important problem is to describe "good" bases of the rings $K(\mathfrak {g})$ as modules over $\mathbb{Z}$ and transition matrices  between them. For example, in the classical case of Lie algebra of type $A(n)$ we have various bases labeled by Young diagrams $\lambda$: Schur polynomials $s_{\lambda}$ (or characters of the irreducible representations), symmetric functions $h_{\lambda}$ and $e_{\lambda}$  (see \cite{Ma} for the details).

We hope also that our result could lead to a better understanding of the algorithms of computing the characters proposed by Serganova and Brundan (see \cite{Brun}, \cite{Serga1}). The investigation of the deformations of the Grothendieck rings and the spectral decompositions of the corresponding analogues of the deformed Calogero-Moser and Macdonald operators \cite{SV} may help to clarify the situation.
%An interesting question also is to find the characters of the Kac modules for the basic %classical Lie superalgebras of the type II.

One can define also the Grothendieck ring $P(\mathfrak {g})$ of projective finite-dimensional $\mathfrak {g}$-modules (cf. Serre \cite{Serre2}). It can be shown that  $P(\mathfrak {g}) \subset K(\mathfrak {g})$
is an ideal in the Grothendieck ring $K(\mathfrak {g}).$ An interesting problem is to describe the structure of $P(\mathfrak {g})$ as a $K(\mathfrak {g})$-module.

\section*{Acknowledgements}
This work has been partially supported by the EPSRC (grant EP/E004008/1),
European Union through the FP6 Marie Curie RTN ENIGMA (Contract
number MRTN-CT-2004-5652) and ESF programme MISGAM.

We are very grateful to Vera Serganova for useful comments, which helped us to improve the paper.
We are grateful also to Andrei Okounkov and to the Mathematics Department of Princeton University for the hospitality in March-April 2007, when most of this work had been completed.

\end{document}